\documentclass{amsart}
\RequirePackage{amsmath}
\RequirePackage{amsthm}
\RequirePackage{amsfonts}
\RequirePackage{amssymb}
\RequirePackage{multirow}
\RequirePackage{multicol}
\RequirePackage{float}
\RequirePackage{enumerate}
\RequirePackage{hyperref}
\RequirePackage{xcolor}

\newtheorem{theorem}{Theorem}[section]
\newtheorem{proposition}[theorem]{Proposition}

\newtheorem{lemma}[theorem]{Lemma}
\theoremstyle{definition}
\newtheorem{definition}[theorem]{Definition}
\newtheorem{example}{Example}

\newcommand{\C}{\mathbb{C}}
\newcommand{\R}{\mathbb{R}}
\newcommand{\D}{\mathbb{D}}
\newcommand{\N}{\mathbb{N}}
\newcommand{\Z}{\mathbb{Z}}
\renewcommand{\H}{\mathbb{H}}
\newcommand{\abs}[1]{\left| #1 \right|}

\title{Characterization of the hyperbolic step of parabolic functions}
\author[M.D. Contreras]{Manuel D. Contreras}
\address{Manuel D. Contreras. Departamento de Matem\'atica Aplicada II and IMUS, Escuela T\'ecnica Superior de Ingenier\'ia, Universidad de Sevilla, Camino de los Descubrimientos, s/n 41092, Sevilla, Spain}
\email{contreras@us.es}
\author[F. J. Cruz-Zamorano]{Francisco J. Cruz-Zamorano}
\address{Francisco J. Cruz-Zamorano. Departamento de Matem\'atica Aplicada II and IMUS, Escuela T\'ecnica Superior de Ingenier\'ia, Universidad de Sevilla, Camino de los Descubrimientos, s/n 41092, Sevilla, Spain}
\email{fcruz4@us.es}
\author[L. Rodr\'iguez-Piazza]{Luis Rodr\'iguez-Piazza}
\address{Luis Rodr\'iguez-Piazza. Departmento de An\'alisis Matem\'atico and IMUS, Facultad de Matem\'aticas, Universidad de Sevilla, Calle Tarfia, s/n 41012 Sevilla, Spain}
\email{piazza@us.es}

\date{\today}
\thanks{This research was supported in part by Ministerio de Innovaci\'on y Ciencia, Spain, project PID2022-136320NB-I00, and Junta de Andaluc\'ia, project P20\_00664. The second author was supported by Ministerio de Universidades, Spain, through the action Ayuda del Programa de Formaci\'on de Profesorado Universitario, reference FPU21/00258; and also by Consejer\'ia de Transformaci\'on Econ\'omica, Industria, Conocimiento y Universidades de la Junta de Andaluc\'ia, Ayuda Predoctoral 2021, reference PREDOC\_00490.}

\begin{document}
\maketitle
\begin{abstract}
A classical problem in Complex Dynamics on hyperbolic domains is to characterize the hyperbolic step of parabolic functions. This topic has been studied by several authors, leading to different results and providing characterizations that depend on the behaviour of the iterates of such function. In this work we provide new characterizations in terms of intrinsic properties of the function. 
\end{abstract}
\tableofcontents
\section{Introduction}
The main goal in Discrete Complex Dynamics is to study the asymptotic properties of the iterated self-composition of a given holomorphic function $g \colon \Omega \to \Omega$, where $\Omega \subset \C_{\infty}$ is a open and connected set in the Riemann sphere. If $\Omega$ is also simply connected, the Riemann Mapping Theorem can be used to show that only three different situations occur, corresponding to the Riemann sphere $\Omega = \C_{\infty}$ (known as the rational case), the complex plane $\Omega = \C$ (entire case) and the unit disc $\Omega = \D$ (hyperbolic case). This work is devoted to the last case.

The historical starting point of the hyperbolic case is due to the following well-known result:
\begin{theorem}[Denjoy-Wolff]
\cite[Theorem 3.2.1]{AbateNewBook}
Let $g \colon \D \to \D$ be a holomorphic map which is neither an elliptic
automorphism nor the identity. Then there exists a point $\tau \in \overline{\D}$ such that $g^n \to \tau$, uniformly on compact sets of $\D$, where
$$g^n = g \circ \overset{(n)}{\cdots} \circ g.$$
Furthermore, if $\tau \in \D$, then $g(\tau) = \tau$. Otherwise, the non-tangential limit of $g$ at $\tau$ is $\tau$. That is,
$$\angle\lim_{z \to \tau}g(z) = \tau.$$
\end{theorem}
With the notation introduced in the previous result, whether $\tau \in \D$ or not, $\tau$ can be seen as a (maybe boundary) fixed point of $g$. The point $\tau$ is called the Denjoy-Wolff point (DW point, for short) of $g$. 

It is common in the literature to use the Denjoy-Wolff Theorem to introduce some terminology. If $\tau \in \D$, then $g$ is usually called elliptic. Otherwise, it is possible to show (see \cite[Corollary 2.5.5]{AbateNewBook}) that the following limit exists:
$$\lambda = \angle\lim_{z \to \tau}\dfrac{g(z)-\tau}{z-\tau} \in (0,1].$$
If $\lambda < 1$, $g$ is called hyperbolic. In the case that $\lambda = 1$, $g$ is said parabolic. 

In this work we will focus in the study of parabolic functions, which can be further classified in two different types. To introduce them, note that the Schwarz-Pick Lemma assures
$$\rho_{\D}(g^{n+2}(z),g^{n+1}(z)) \leq \rho_{\D}(g^{n+1}(z),g^n(z)), \quad z \in \D, \, n \in \N,$$
where $\rho_{\D}$ is the pseudo-hyperbolic distance in $\D$. That is, the sequence of distances $\{\rho_{\D}(g^{n+1}(z),g^n(z))\}$ is non-negative and non-increasing. Then, it must be convergent. Indeed (see \cite[Corollary 4.6.9.(i)]{AbateNewBook}), if there exists $z \in \D$ such that $\rho_{\D}(g^{n+1}(z),g^n(z)) \to 0$, then $\rho_{\D}(g^{n+1}(w),g^n(w)) \to 0$ for every $w \in \D$. As this property only depends on $g$, the following concept can be defined:
\begin{definition}
\label{def:hypstep}
If a holomorphic function $g \colon \D \to \D$ is such that there exists $z \in \D$ with $\rho_{\D}(g^{n+1}(z),g^n(z)) \to 0$, then $g$ is said to be of zero hyperbolic step (0HS, for short). Otherwise, is it said of positive hyperbolic step (PHS).
\end{definition}
The hyperbolic step of a function plays an important role in its dynamical properties. For example, elliptic functions are always of 0HS, unless they are automorphisms (for automorphisms, $\rho_{\D}(g^{n+1}(z),g^n(z))$ is a constant sequence). On the other hand, hyperbolic functions are always of PHS. However, the situation for parabolic functions is different. That is, it is possible to give explicit examples (see Section \ref{sec:examples}) of some parabolic functions with 0HS and some others with PHS (even if those functions are not automorphisms).

A classical problem in the theory of iteration of holomorphic self-maps is to achieve a complete characterization of the hyperbolic step of parabolic functions, and much work has been done in this sense. Some authors have given a complete characterization in terms of different concepts that are related to parabolic functions. One of them are orbits. For example, Cowen showed that a parabolic function is of PHS if and only if its orbits are interpolating sequences \cite[Proposition 4.9]{Cowen83} and Pommerenke provided a different characterization in terms of the limit behaviour of an orbit (see the results from \cite{PommerenkeHalfPlane} that are cited in Section \ref{sec:preliminaries}). Other authors have used models to characterize the hyperbolic step \cite[Section 3]{Cowen81} or geometric properties of the Koenigs function \cite[Proposition 3.3]{Abel}. However, there is not a complete characterization in terms of intrisic properties of the function. Both in the paper of Bourdon and Shapiro \cite{BourdonShapiro} and in the paper of Contreras, Díaz-Madrigal and Pommerenke \cite{CDP}, there are different criteria to classify parabolic functions assuming some extra regularity of the function at the Denjoy-Wolff point. See the book of Abate \cite[Section 4.7]{AbateNewBook} for a selection of such results. It is worth mentioning that an inner function is of 0HS if and only if the corresponding boundary map is ergodic, as noted in \cite[Section 3]{Doering}. Indeed, there are also some interesting results for inner functions, as the ones due to Aaronson (see the different results from \cite{Aaronson} that are cited through the document). There also some references on this topic in the context of continuous semigroups, see the recent monograph \cite{Contreras}.

The purpose of this paper is to provide intrinsic characterizations of the hyperbolic step for a given parabolic function. To do so, first note that working with $\D$ as a domain is not crucial, as everything can be translated to any other hyperbolic domain. In this work we will use holomorphic self-maps of the upper halfplane $\H = \{z \in \C \mid \text{Im}(z) > 0\}$. Indeed, there is a correspondence between every function $g \colon \D \to \D$ and the function $f \colon \H \to \H$ given by $f = S^{-1} \circ g \circ S$, where $S \colon \D \to \H$ is any bijective Möbius transformation. These new functions have a useful description, as it is known (see \cite[Theorem 6.2.1]{Aaronson}) that every holomorphic function $f \colon \H \to \H$ can be uniquely written as
$$f(z) = \alpha z + \beta + \int_{\R}\dfrac{1+tz}{t-z}d\mu(t), \quad z \in \H,$$
where $\alpha \geq 0$, $\beta \in \R$ and $\mu$ is a positive finite measure on $\R$. In this context, the Denjoy-Wolff theorem can be obviously rewritten. For example, if $f$ is not the identity map and $\alpha \geq 1$, $f$ can not have fixed points in $\H$ and $\text{Im}(f^n)$ is a non-decreasing sequence, therefore the corresponding DW point is infinity. Indeed, one has (see \cite[Chapter 5, Lemma 2]{Doering}) $f(it)/(it) \to \alpha$ as $t \to +\infty$. This assures that $\alpha \geq 1$ if and only if the DW point of $f$ is infinity. Actually, in this situation, $f$ is parabolic if and only if $\alpha = 1$, so the following is derived:
\begin{theorem}
\label{teo:expresion}
Every parabolic function $f \colon \H \to \H$ whose DW point is infinity can be uniquely written as
\begin{equation}
\label{eq:rep}
f(z) = z + \beta + \int_{\R}\dfrac{1+tz}{t-z}d\mu(t), \quad z \in \H,
\end{equation}
where $\beta \in \R$ and $\mu$ is a positive finite measure on $\R$. Conversely, any function of the form \eqref{eq:rep} is parabolic and infinity is its DW point whenever $\beta$ and $\mu$ are not simultaneously zero.
\end{theorem}
By means of a conjugation, one can always suppose that a parabolic map has its DW point at infinity, so the previous result can be used to describe the behaviour of all the parabolic functions. In fact, every dynamical property of a parabolic function $f$ can only depend on its associated parameters $\beta$ and $\mu$. In particular, in this work we give several characterizations of the hyperbolic step of $f$ in terms of these quantities under different scenarios.

In Section \ref{sec:t1}, we give a complete description of the hyperbolic step assuming $t \in L^1(\mu)$, that is,
$$\int_{\R}\abs{t}d\mu(t) < +\infty.$$
Indeed, some other integrability conditions show to be crucial for the hyperbolic step of the corresponding function. In that section, the main result can be summarized as follows:
\begin{theorem}
\label{teo:t1end}
Let $f \colon \H \to \H$ be a parabolic function whose DW point is infinity, given by Theorem \ref{teo:expresion}. Suppose that $t \in L^1(\mu)$.
\begin{enumerate}[\hspace{0.5cm}\normalfont(a)]
\item If $\beta = \displaystyle\int_{\R}td\mu(t)$, then $f$ is of 0HS.
\item If $\beta > \displaystyle\int_{\R}td\mu(t)$, then $f$ is of PHS if and only if $\displaystyle\int_{(0,+\infty)}t^2d\mu(t) < +\infty$.
\item If $\beta < \displaystyle\int_{\R}td\mu(t)$, then $f$ is of PHS if and only if $\displaystyle\int_{(-\infty,0)}t^2d\mu(t) < +\infty$.
\end{enumerate}
Furthermore, given $z \in \H$, $f$ is of 0HS if and only if $\textnormal{Im}(f^n(z)) \to +\infty$.
\end{theorem}
This result generalizes a result from Aaronson (see Theorem \ref{teo:AaronsonCompacto}), who provided the characterization assuming that $\mu$ is of compact support.

In Section \ref{sec:symmetric} we also work on an idea that was previously studied by Aaronson, that is, functions corresponding to symmetric measures. In this context, we prove the following extension of another theorem from Aaronson (see \cite[Theorem 6.4.5]{Aaronson}):
\begin{theorem}
\label{teo:symmetric_conc}
Let $f \colon \H \to \H$ be a parabolic function whose DW point is infinity, given by Theorem \ref{teo:expresion}. Suppose that $\mu$ is symmetric.
\begin{enumerate}[\hspace{0.5cm}\normalfont(a)]
\item If $\displaystyle\int_{\R}t^2d\mu(t) < +\infty$, then $f$ is of 0HS if and only if $\beta = 0$.
\item If $\displaystyle\int_{\R}t^2d\mu(t) = +\infty$, then $f$ is of 0HS for every $\beta \in \R$.
\end{enumerate}
Furthermore, given $z \in \H$, $f$ is of 0HS if and only if $\textnormal{Im}(f^n(z)) \to +\infty$.
\end{theorem}

In Section \ref{sec:perturb} we start by characterizing the hyperbolic step of parabolic functions whose corresponding measure is supported on a half-line (see Proposition \ref{prop:weak}). Indeed, two different proofs of this are given. The first of them is shown in Subsection \ref{subsec:firstproof}, and it uses the concept of Koenigs function. The three authors want to thank Pavel Gumenyuk for some fruitful discussions with him, which eventually led to this proof. In Subsection \ref{subsec:secondproof} we give a different proof, whose ideas will be later improved in Subsection \ref{subsec:main} to show that the following more general situation holds:
\begin{theorem}
\label{teo:perturb}
Let $f \colon \H \to \H$ be a parabolic function whose DW point is infinity, given by Theorem \ref{teo:expresion}. Suppose that $\mu$ is such that
$$\int_{(-\infty,0)}\abs{t}d\mu(t) = \infty, \quad \int_{(0,+\infty)}t^2d\mu(t) < \infty.$$
Then $f$ is of PHS for any $\beta \in \R$.
\end{theorem}
Note that this result is decidedly different from the previous two characterizations: in this one there exists no $\beta \in \R$ such that the corresponding $f$ is of 0HS.

To end, Section \ref{sec:examples} is devoted to explicitly state several examples in order to clarify our results. We begin with Section \ref{sec:preliminaries}, where some previous results on the hyperbolic step of parabolic functions are revisited.

The authors thank the referee for his/her carefully reading of the original manuscript.
\section{Preliminaries}
\label{sec:preliminaries}
For the rest of the document, $f \colon \H \to \H$ denotes a parabolic function whose DW point is infinity. Similarly, $\{z_n\} \subset \H$ will denote the orbit of $z_0 \in \H$ by $f$, where
$$z_n = f^n(z_0), \quad x_n = \text{Re}(z_n), \quad y_n = \text{Im}(z_n), \quad n \in \N.$$
Note that $y_n$ is an increasing sequence of positive numbers, and so its limit exists, although it might be infinite.

With this notation, the following holds:
\begin{theorem}
\cite{PommerenkeHalfPlane}
\label{teo:Pomm}
The following limit
$$b = \lim_{n \to \infty}\dfrac{x_{n+1}-x_n}{y_n} \in \R$$
exists. Furthermore, $b = 0$ if and only if $f$ is of 0HS. Also, $z_{n+1}/z_n \to 1$ (see \cite[Eq. (3.16)]{PommerenkeHalfPlane}), and $y_{n+1}/y_n \to 1$ (see \cite[Eq. (3.17)]{PommerenkeHalfPlane}).
\end{theorem}
\begin{theorem}
\label{teo:CDP}
\begin{enumerate}[\normalfont(a)]
\item \cite[Theorem 4.1]{CDP} If $f$ is of PHS, then the following are equivalent:
\subitem(i) $\displaystyle\angle\lim_{z \to \infty} (f(z)-z) \in \C$.
\subitem(ii) $\displaystyle\lim_{n \to \infty}y_n < + \infty$.
\item \cite[Proposition 3.3]{CDP} If $f$ is of 0HS, then $\lim_{n \to \infty}y_n = + \infty$.
\end{enumerate}
\end{theorem}
\begin{theorem}
\label{teo:convtang}
\cite[Remark 1]{PommerenkeHalfPlane}
If $f$ is of PHS, then $z_n$ converges tangentially to infinity. That is, $\arg(z_n) \to 0$ or $\arg(z_n) \to \pi$.
\end{theorem}
\section{Characterization when $t \in L^1(\mu)$}
\label{sec:t1}
Let $f$ be given by the expression of Theorem \ref{teo:expresion}, and suppose that $\mu$ is such that
$$\int_{\R} \abs{t} d\mu(t) < \infty.$$
Then, $f$ can be rewritten as
\begin{equation}
\label{eq:ft1}
f(z) = z + \tilde{\beta} + \int_{\R}\dfrac{1+t^2}{t-z}d\mu(t), \quad z \in \H,
\end{equation}
where
$$\tilde{\beta} = \beta - \int_{\R}td\mu(t).$$

Furthermore, $f$ can also be rewritten as
$$f(z) = z + \tilde{\beta} + \int_{\R}\dfrac{d\nu(t)}{t-z}, \quad z \in \H,$$
where $d\nu(t) = (1+t^2)d\mu(t)$, and so $\nu$ is a finite measure if and only if
$$\int_{\R} t^2 d\mu(t) < \infty.$$

Indeed, this situation can be deduced directly from the properties of $f$ as follows:
\begin{proposition}
\label{prop:t2charac}
Let $f \colon \H \to \H$ be a parabolic function whose Denjoy-Wolff's point is infinity, given by Theorem \ref{teo:expresion}. The following statements are equivalent:
\begin{enumerate}[\hspace{0.5cm}\normalfont(i)]
\item $\displaystyle\int_{\R}t^2d\mu(t) < \infty$.
\item There exists $\beta \in \R$ such that $\displaystyle\angle\lim_{z \to \infty}z(f(z)-z) \in (-\infty,0)$.
\end{enumerate}
Moreover, if such $\beta$ exists, then it must be
$$\beta = \int_{\R}td\mu(t), \quad \angle\lim_{z \to \infty}z(f(z)-z) = -\int_{\R}(1+t^2)d\mu(t).$$
\begin{proof}
Suppose
$$\int_{\R}t^2d\mu(t) < \infty.$$
Then, as said before, $f$ can be written as
$$f(z) = z + \tilde{\beta} + \int_{\R}\dfrac{d\nu(t)}{t-z}, \quad z \in \H,$$
where $\nu$ is a positive finite measure on $\R$.

Note that if
$$z = x+iy \in S = \{w \in \H : \abs{\text{Re}(w)} \leq C\text{Im}(w)\}$$
for some $C > 0$, then
$$\abs{\dfrac{z}{t-z}}^2 = \dfrac{x^2+y^2}{(t-x)^2+y^2} \leq \dfrac{x^2+y^2}{y^2} \leq C^2+1.$$
Thus, it is a mere consequence of Lebesgue's Dominated Convergence Theorem that
$$\angle\lim_{z \to \infty}z\int_{\R}\dfrac{d\nu(t)}{t-z} = -\nu(\R).$$
Consequently,
$$\angle\lim_{z \to \infty}z(f(z)-z) \in (-\infty,0) \iff \tilde{\beta} = 0,$$
that is,
$$\beta = \int_{\R}td\mu(t).$$
In that case,
$$\angle\lim_{z \to \infty}z(f(z)-z) = -\nu(\R) = -\int_{\R}(1+t^2)d\mu(t).$$

On the other hand, suppose that
$$\angle\lim_{z \to \infty}z(f(z)-z) \in (-\infty,0).$$
In particular,
$$\lim_{y \to +\infty} iy(f(iy)-iy) < \infty,$$
and so
$$\lim_{y \to +\infty} \text{Re}\left(iy(f(iy)-iy)\right) = \lim_{y \to + \infty} -\int_{\R}\dfrac{(1+t^2)y^2}{t^2+y^2}d\mu(t) > - \infty.$$
But, using Fatou's Lemma, one has
$$\lim_{y \to + \infty} \int_{\R}\dfrac{(1+t^2)y^2}{t^2+y^2}d\mu(t) \geq \int_{\R}(1+t^2)d\mu(t).$$
From there, we deduce that
$$\int_{\R}t^2d\mu(t) < \infty.$$
\end{proof}
\end{proposition}

In this section, our interest is to characterize the hyperbolic step of functions described by \eqref{eq:ft1}. This question has been partially studied, and the following is known:
\begin{theorem}
\label{teo:AaronsonCompacto}
\cite[Theorem 6.4.1]{Aaronson}
Let $f \colon \H \to \H$ be an inner function given by
$$f(z) = z + \tilde{\beta} + \int_{\R}\dfrac{d\nu(t)}{t-z}, \quad z \in \H,$$
where $\tilde{\beta} \in \R$ and $\nu$ is a positive and finite measure of compact support. Then $f$ is ergodic if and only if $\tilde{\beta}= 0$.
\end{theorem}
For inner functions (those satisfying $\angle\lim_{z \to x}\text{Im}(f(z)) = 0$ for almost all $x\in \R$) the properties of ergodicity and being of 0HS are equivalent (see \cite[Theorem 3.1]{Doering}). Then, the previous theorem is a partial answer for our characterization.

To give a complete description, we need some previous lemmas. In order to introduce them, let us define
$$p(z) = f(z) - z - \tilde{\beta} = \int_{\R}\dfrac{1+t^2}{t-z}d\mu(t), \quad z \in \H.$$
Then, the following holds:
\begin{lemma}
\label{lemma:anglep0}
If $t \in L^1(\mu)$, then  $\angle\lim_{z \to \infty}p(z) = 0$.
\begin{proof}
Using Lindelöf's Theorem, it suffices to prove that
$$\lim_{y \to + \infty}p(iy) = 0.$$
To do so, notice that
$$p(iy) = \int_{\R}\dfrac{1+t^2}{t-iy}d\mu(t),$$
where
$$\abs{\dfrac{1+t^2}{t-iy}} = \dfrac{1+t^2}{\sqrt{t^2+y^2}} \leq \dfrac{1+t^2}{\sqrt{t^2+1}}, \quad y \geq 1.$$
As $t \in L^1(\mu)$, then
$$\int_{\R}\dfrac{1+t^2}{\sqrt{t^2+1}}d\mu(t)<\infty.$$
So the result follows from Lebesgue's Dominated Convergence Theorem.
\end{proof}
\end{lemma}
By using the ideas in \cite{CDP}, we have the following result:
\begin{proposition}
\label{prop:phsiffconv}
Suppose $t \in L^1(\mu)$. Then, $f$ is of PHS if and only if $\lim_n y_n < +\infty$.
\begin{proof}
The result follows from Theorem \ref{teo:CDP} by using Lemma \ref{lemma:anglep0}.
\end{proof}
\end{proposition}
In order to characterize the hyperbolic step when $\tilde{\beta} = 0$, we need the following result:
\begin{lemma}
\label{lemma:phs}
Suppose that $t \in L^1(\mu)$ and $f$ is of PHS, and define $q(z) = f(z) - z$ for $z \in \H$. Then, 
\begin{enumerate}[\hspace{0.5cm}\normalfont(a)]
\item $\displaystyle\lim_{n \to \infty}q(z_n) = \angle\lim_{z \to \infty}q(z) \in \R \setminus \{0\}$.
\item $p(z_n) \to 0$.
\end{enumerate}
\begin{proof}
Suppose $f$ is of PHS. Then, by Proposition \ref{prop:phsiffconv}, there must exists $L > 0$ such that $y_n \to L$. In that case,
$$\text{Im}(q(z_n)) = y_{n+1}-y_n = y_n\left(\dfrac{y_{n+1}}{y_n}-1\right) \to 0,$$
as $y_n \to L$ and $y_{n+1}/y_n \to 1$, by using Theorem \ref{teo:Pomm}. But $q(z) = \tilde{\beta} + p(z)$, so $\text{Im}(q(z_n)) = \text{Im}(p(z_n))$, and then $\text{Im}(p(z_n)) \to 0$.

Concerning the real part, in \cite[Lemma 3.5]{CDP} it is shown that there exists $C > 0$ such that
$$\abs{q(z)-q(z_n)} \leq C\text{Im}(q(z_n)), \quad z \in I_n,$$
where
$$I_n = \{tz_n + (1-t)z_{n+1} : 0 \leq t \leq 1\}.$$
So,
$$\sup_{z \in I_n} \abs{q(z)-q(z_n)} \leq C \text{Im}(q(z_n)) \to 0.$$
Then, define $\Gamma = \bigcup_{n \in \N}I_n$ and notice that
\begin{align*}
\lim_{\substack{z \to \infty \\ z \in \Gamma}}q(z) & = \lim_{n \to \infty}q(z_n) = \lim_{n \to \infty} z_{n+1}-z_n \\
& = \lim_{n \to \infty}y_n\left(\dfrac{x_{n+1}-x_n}{y_n}+i\dfrac{y_{n+1}-y_n}{y_n}\right) \\
& = L(b+i0) = Lb,
\end{align*}
where Theorem \ref{teo:Pomm} is used and $L > 0$, $b \neq 0$.

Using Lindelöf's Theorem, it follows that $\angle\lim_{z \to \infty}q(z) = bL \in \R \setminus \{0\}$. Actually, using Lemma \ref{lemma:anglep0}, $\tilde{\beta} = bL$. Furthermore,
$$bL = \lim_{n \to \infty}\dfrac{x_{n+1}-x_n}{y_n}\lim_{n \to \infty}y_n = \lim_{n \to \infty} x_{n+1}-x_n.$$
So,
$$x_{n+1}-x_n = \tilde{\beta}+\text{Re}(p(z_n)) \to bL = \tilde{\beta},$$
and then, $\text{Re}(p(z_n)) \to 0$.
\end{proof}
\end{lemma}
With the former ideas, the following can be shown:
\begin{proposition}
\label{prop:beta0}
Suppose $t \in L^1(\mu)$. If $\tilde{\beta} = 0$, then $f$ is of 0HS.
\begin{proof}
If $\tilde{\beta} = 0$, then $q(z) = f(z) - z = p(z)$. In that case, if $f$ is supposed to be of PHS, then Lemma \ref{lemma:anglep0} provides a contradiction.
\end{proof}
\end{proposition}
Now, it remains to characterize the case where $\tilde{\beta} \neq 0$, which is more subtle. The hyperbolic step of $f$ when $\tilde{\beta} \neq 0$ depends on some integrability conditions with respect to $\mu$. For example, one can prove the following:
\begin{theorem}
\label{teo:t2case}
Let $f \colon \H \to \H$ be a parabolic function given by
$$f(z) = z + \tilde{\beta} + \int_{\R}\dfrac{d\nu(t)}{t-z}, \quad z \in \H,$$
where $\tilde{\beta} \in \R$ and $\nu$ is a positive and finite measure on $\R$. Then $f$ is of 0HS if and only if $\tilde{\beta} = 0$.
\begin{proof}
By Proposition \ref{prop:beta0}, it remains to prove that $f$ is of PHS when $\tilde{\beta} \neq 0$. By Proposition \ref{prop:phsiffconv}, it is enough to prove that $\lim_n y_n < \infty$. To do so, define
$$\varepsilon_n = \dfrac{y_{n+1}-y_n}{y_n} > 0,$$
and note that $y_{n+1} = y_n + (y_{n+1} - y_n) = y_n(1+\varepsilon_n)$. Then,
$$y_n = y_0 \prod_{k = 0}^{n-1}(1+\varepsilon_k).$$
Thus,
$$\lim_{n \to \infty}y_n < \infty \iff \prod_{n = 0}^{\infty}(1+\varepsilon_n) < \infty \iff \sum_{n = 0}^{\infty}\varepsilon_n < \infty.$$

On the other hand,
$$\varepsilon_n = \int_{\R}\dfrac{d\nu(t)}{(t-x_n)^2+y_n^2}.$$
In that case, it is enough to prove that there exists $C > 0$ such that
$$\sum_{n = 0}^{\infty}\dfrac{1}{(t-x_n)^2+y_n^2} \leq C < \infty.$$

To do so, define
$$p(z) = \int_{\R}\dfrac{d\nu(t)}{t-z}, \quad z \in \H.$$
Note that
$$\abs{\dfrac{1}{t-z}} \leq \dfrac{1}{\text{Im}(z)}, \quad z \in \H.$$
Consequently, as $z_n \to \infty$, it is a mere application of Lebesgue's Dominated Convergence Theorem that $p(z_n) \to 0$.

In that case, suppose that $\tilde{\beta} > 0$ (similar ideas apply if $\tilde{\beta} < 0$) and assume that
$$\abs{p(z_n)} \leq \dfrac{\tilde{\beta}}{2}, \quad n \in \N.$$
Then,
$$\dfrac{\tilde{\beta}}{2} \leq x_{n+1}-x_n = \tilde{\beta} + \text{Re}(p(z_n)) \leq \dfrac{3\tilde{\beta}}{2}, \quad n \in \N.$$
In particular, $x_n$ is an increasing sequence. From this, it can be deduce that
$$x_n+\dfrac{\tilde{\beta}}{2}(m-n) \leq x_m \leq x_n +\dfrac{3\tilde{\beta}}{2}(m-n), \quad m \geq n \in \N.$$

So, if $t \leq x_0$, then
$$x_n-t = (x_n-x_0)+(x_0-t) \geq x_n - x_0 \geq \dfrac{\tilde{\beta}}{2}n.$$
In that case,
$$\sum_{n \in \N}\dfrac{1}{(t-x_n)^2+y_n^2} \leq \sum_{n \in \N}\dfrac{1}{n^2\tilde{\beta}^2/4+y_0^2} < \infty.$$

Otherwise, if $t > x_0$, it is possible to (uniquely) find $N \in \N$ such that $t \in [x_N,x_{N+1})$. Then, for $n \in \N$ with $n \leq N$,
$$t-x_n = (t - x_N) + (x_N-x_n) \geq x_N-x_n \geq (N-n)\dfrac{\tilde{\beta}}{2}.$$
On the other hand, if $n > N$, then
$$x_n - t = (x_n-x_{N+1})+(x_{N+1}-t) \geq x_n-x_{N+1} \geq (n-N-1)\dfrac{\tilde{\beta}}{2}.$$
Joining both estimates,
\begin{align*}
\sum_{n \in \N}\dfrac{1}{(t-x_n)^2+y_n^2} & \leq \sum_{n \in \N}\dfrac{1}{(t-x_n)^2+y_0^2} \\
& = \sum_{n = 0}^N\dfrac{1}{(t-x_n)^2+y_0^2} +\sum_{n = N+1}^{\infty}\dfrac{1}{(t-x_n)^2+y_0^2} \\
& \leq \sum_{n = 0}^N\dfrac{1}{(N-n)^2\tilde{\beta}^2/4+y_0^2} +\sum_{n = N+1}^{\infty}\dfrac{1}{(n-N-1)^2\tilde{\beta}^2/4+y_0^2} \\
& \leq 2\sum_{n \in \N}\dfrac{1}{n^2\tilde{\beta}^2/4+y_0^2} < \infty.
\end{align*}
\end{proof}
\end{theorem}

The previous theorem is a complete characterization of the hyperbolic step of $f$ when $t \in L^2(\mu)$, so that it enhances Aaronson's Theorem \ref{teo:AaronsonCompacto}. Using this integrability condition, one can prove the following:
\begin{proposition}
\label{prop:necessary}
Suppose $t \in L^1(\mu)$.
\begin{enumerate}[\normalfont(i)]
\item If $\tilde{\beta} > 0$ and $f$ is of PHS, then $\displaystyle\int_{(0,+\infty)}t^2d\mu(t) < \infty$.
\item If $\tilde{\beta} < 0$ and $f$ is of PHS, then $\displaystyle\int_{(-\infty,0)}t^2d\mu(t) < \infty$.
\end{enumerate}
\begin{proof}
As the ideas are similar, we will only prove (\textit{i}). If $f$ is of PHS, then Proposition \ref{prop:phsiffconv} shows that the sequence $\{y_n\}$ is convergent, then bounded. But
$$y_{n+1} = y_n\left(1+\dfrac{y_{n+1}-y_n}{y_n}\right) = y_0\prod_{k=0}^n\left(1+\dfrac{y_{k+1}-y_k}{y_k}\right).$$
Then
$$y_n \text{ is bounded} \iff \prod_{n=0}^{\infty}\left(1+\dfrac{y_{n+1}-y_n}{y_n}\right) < +\infty \iff \sum_{n = 0}^{\infty}\dfrac{y_{n+1}-y_n}{y_n} < +\infty,$$
where
$$\sum_{n = 0}^{\infty}\dfrac{y_{n+1}-y_n}{y_n} = \int_{\R}\left(\sum_{n = 0}^{\infty}\dfrac{1+t^2}{(t-x_n)^2+y_n^2}\right)d\mu(t).$$

On the other hand, $p(z_n) \to 0$, as shown in Lemma \ref{lemma:phs}. In that case, $x_{n+1}-x_n \to \tilde{\beta} > 0$. Then, let $N \in \N$ be such that
$$\dfrac{\tilde{\beta}}{2} \leq x_{n+1}-x_n \leq \dfrac{3\tilde{\beta}}{2}, \quad n \geq N.$$
Thus, for every $t > x_N$ one can find $k = k(t) \geq N$ such that $\abs{t-x_k} \leq \tilde{\beta}$. In that case,
$$\sum_{n \in \N}\dfrac{1+t^2}{(t-x_n)^2+y_n^2} \geq \dfrac{1+t^2}{\tilde{\beta}^2+y_k^2}, \quad t > x_N.$$

Furthermore, let $L > 0$ be such that $y_n \to L$, as Proposition \ref{prop:phsiffconv} assures. Then,
$$\dfrac{1}{\tilde{\beta}^2+y_k^2} \to \dfrac{1}{\tilde{\beta}^2+L^2} = C.$$
In particular,
$$\sum_{n \in \N}\dfrac{1+t^2}{(t-x_n)^2+y_n^2} \geq C(1+t^2), \quad t > x_N.$$
And so,
$$+\infty > \int_{\R}\left(\sum_{n = 0}^{\infty}\dfrac{1+t^2}{(t-x_n)^2+y_n^2}\right)d\mu(t) \geq \int_{(x_N,+\infty)}C(1+t^2)d\mu(t).$$
The result easily follows from these ideas.
\end{proof}
\end{proposition}

By joining the ideas of Proposition \ref{prop:beta0}, Theorem \ref{teo:t2case} and Proposition \ref{prop:necessary}, it only remains to characterize the hyperbolic step of $f$ in two similar scenarios. For $\tilde{\beta} > 0$ under the hypotheses
$$\int_{\R}\abs{t}d\mu(t) < +\infty, \quad \int_{(-\infty,0)}t^2d\mu(t) = +\infty, \quad \int_{(0,+\infty)}t^2d\mu(t) < +\infty,$$
and for $\tilde{\beta} < 0$ under the hypotheses
$$\int_{\R}\abs{t}d\mu(t) < +\infty, \quad \int_{(-\infty,0)}t^2d\mu(t) < +\infty, \quad \int_{(0,+\infty)}t^2d\mu(t) = +\infty.$$
Until the end of the section, we will focus on proving that, under both of these settings, $f$ is of PHS. Indeed, it is enough to do so for the first of them. Note that we can decompose $\mu$ by $\mu = \mu_1+\mu_2$, where
\begin{equation}
\label{eq:decomp}
\mu_1(A) = \mu(A \cap (-\infty,0]), \quad \mu_2(A) = \mu(A \cap (0,+\infty)).
\end{equation}
In particular,
$$\int_{\R}t^2d\mu_2(t) < +\infty.$$
Similarly, decompose $p$ by $p = p_1 + p_2$ where
$$p_j(z) = \int_{\R}\dfrac{1+t^2}{t-z}d\mu_j(t), \quad z \in \H, \, j \in \{1,2\}.$$
We begin by showing the following result:
\begin{lemma}
\label{lemma:lastcase}
Under the previous hypotheses, there exist $a,b > 0$ large enough so that
$$\Omega = \{x+iy \in \H : x > a, y > b\}$$
is invariant by $f$, that is, $f(\Omega) \subset \Omega$. Furthermore, if $z_0 \in \Omega$, then $p(z_n) \to 0$.
\begin{proof}
First, note that if $z = x+iy \in \H$ then
$$\abs{\dfrac{1+t^2}{t-z}} \leq \dfrac{1}{y}(1+t^2) \leq 1+t^2, \quad y \geq 1.$$
As the last bound is integrable with respect to $\mu_2$, it is a mere consequence of Lebesgue's Dominated Convergence Theorem that
$$\lim_{\substack{z \to \infty \\ y > 1}}p_2(z) = 0.$$

On the other hand, 
$$\abs{p_1(z)} \leq \int_{(-\infty,0]}\dfrac{1+t^2}{\abs{t-z}}d\mu(t) \leq \int_{(-\infty,0]}(1+\abs{t})\abs{\dfrac{t-i}{t-z}}d\mu(t).$$
It is easy to see that if $y > 1$, $x >0$ and $t \leq 0$, then
$$\abs{\dfrac{t-i}{t-z}}^2 = \dfrac{1+t^2}{(t-x)^2+y^2} \leq \dfrac{1+t^2}{t^2+y^2} \leq 1.$$
Also,
$$\lim_{t \to - \infty}\abs{\dfrac{t-i}{t-z}} = 1,$$
so
$$\sup_{t \leq 0} \abs{\dfrac{t-i}{t-z}} = 1.$$
In that case, Lebesgue's Dominated Convergence Theorem shows that
$$\lim_{\substack{z \to \infty \\ x > 0, \, y > 1}} p_1(z) = 0.$$

Applying the previous ideas, we can choose $a,b > 0$ such that
$$\abs{p_1(z)} \leq \dfrac{\tilde{\beta}}{4}, \quad \abs{p_2(z)} \leq \dfrac{\tilde{\beta}}{4}, \quad z \in \Omega,$$
where
$$\Omega = \{x+iy \in \H : x > a, y > b\}.$$
Thus, if $z \in \Omega$, then
$$\text{Re}(f(z)) - x = \tilde{\beta} + \text{Re}(p_1(z)) + \text{Re}(p_2(z)) \geq \dfrac{\tilde{\beta}}{2} > 0.$$
This shows that $f(\Omega) \subset \Omega$, where we have used that $\text{Im}(f(z)) > \text{Im}(z)$ for all $z \in \H$. So, if $z_0 \in \Omega$, then $z_n \in \Omega$ for all $n \in \N$, and then
$$\lim_{n \to \infty}p(z_n) = \lim_{\substack{z \to \infty \\ z \in \Omega}}p(z) = 0,$$
as shown before.
\end{proof}
\end{lemma}
Remember that, by Proposition \ref{prop:phsiffconv}, to see that $f$ is of PHS it suffices to show that $y_n$ is bounded. But, as we have used several times before, one can see that
$$y_{n+1} = y_n\left(1+\dfrac{y_{n+1}-y_n}{y_n}\right) = y_0 \prod_{k = 0}^n\left(1+\dfrac{y_{k+1}-y_k}{y_k}\right).$$
Then, 
\begin{align*}
\sup_{n \in \N} y_n < \infty & \iff \prod_{n = 0}^{\infty}\left(1+\dfrac{y_{n+1}-y_n}{y_n}\right) < \infty \\
& \iff \sum_{n = 0}^{\infty}\dfrac{y_{n+1}-y_n}{y_n} < \infty \\
& \iff \sum_{n = 0}^{\infty}\int_{\R}\dfrac{1+t^2}{(t-x_n)^2+y_n^2}d\mu(t) < \infty.
\end{align*}

Also, by using Lemma \ref{lemma:lastcase}, we can suppose that $z_0$ is such that
$$\abs{p(z_n)} \leq \dfrac{\tilde{\beta}}{2}, \quad n \in \N.$$
So,
$$\dfrac{\tilde{\beta}}{2} \leq x_{n+1}-x_n \leq \dfrac{3\tilde{\beta}}{2}, \text{ and, } x_0+n\dfrac{\tilde{\beta}}{2} \leq x_n \leq x_0+n\dfrac{3\tilde{\beta}}{2}.$$
Indeed, we can suppose that
$$x_0 > \dfrac{\tilde{\beta}}{2}.$$

With this ideas, and using the decomposition for $\mu$ given in \eqref{eq:decomp}, we can prove the following:
\begin{lemma}
\label{lemma:eps1}
Under the previous hypotheses,
$$\sum_{n = 0}^{\infty}\int_{\R}\dfrac{1+t^2}{(t-x_n)^2+y_n^2}d\mu_j(t) < \infty, \quad j \in \{1,2\}.$$
\begin{proof}
For $j = 2$, the ideas are the same than the ones given in the proof of Theorem \ref{teo:t2case}. For $j = 1$, the proof goes as follows:

Note that, if $t \leq 0$,
$$\dfrac{1+t^2}{(t-x_n)^2+y_n^2} \leq \dfrac{1+t^2}{(t-x_n)^2} = \dfrac{1+t^2}{(\abs{t}+x_n)^2} \leq \dfrac{1+t^2}{(\abs{t}+x_0+n\tilde{\beta}/2)^2}.$$
In that case,
\begin{align*}
\sum_{n \in \N}\dfrac{1+t^2}{(t-x_n)^2+y_n^2} & \leq \sum_{n \in \N}\dfrac{1+t^2}{(\abs{t}+x_0+n\tilde{\beta}/2)^2} \\
& = \dfrac{4(1+t^2)}{\tilde{\beta}^2}\sum_{n \in \N}\dfrac{1}{(2(\abs{t}+x_0)/\tilde{\beta}+n)^2} \\
& \leq \dfrac{4(1+t^2)}{\tilde{\beta}^2}\int_{-1}^{\infty}\dfrac{1}{(2(\abs{t}+x_0)/\tilde{\beta}+s)^2}ds.
\end{align*}

Now, call $\alpha = 2(\abs{t}+x_0)/\tilde{\beta}$. Note that $\alpha \geq 2x_0/\tilde{\beta} > 1$, as $x_0$ is chosen so that $x_0 > \tilde{\beta}/2$. With this notation, we have
$$\sum_{n \in \N}\dfrac{1+t^2}{(t-x_n)^2+y_n^2} \leq \dfrac{4(1+t^2)}{\tilde{\beta}^2}\int_{-1}^{\infty}\dfrac{1}{(\alpha+s)^2}ds,$$
where $\alpha + s > 0$ for $s > -1$. So,
\begin{align}
\sum_{n \in \N}\dfrac{1+t^2}{(t-x_n)^2+y_n^2} & \leq \dfrac{4(1+t^2)}{\tilde{\beta}^2}\left[-\dfrac{1}{\alpha+s}\right]^{\infty}_{-1} \label{eq:ineqsum}  \\
& = \dfrac{4(1+t^2)}{\tilde{\beta}^2}\dfrac{1}{\alpha-1} = \dfrac{2}{\tilde{\beta}}\dfrac{1+t^2}{\abs{t}+x_0-\tilde{\beta}/2}. \nonumber
\end{align}
Once more, note that
$$x_0 > \dfrac{\tilde{\beta}}{2}, \text{ and, } \abs{t}+x_0-\dfrac{\tilde{\beta}}{2} > 0, \quad t \leq 0.$$

To sum up, remember that $\mu_1$ is supported on $(-\infty,0]$, and then
\begin{align*}
\sum_{n = 0}^{\infty}\int_{\R}\dfrac{1+t^2}{(t-x_n)^2+y_n^2}d\mu_1(t) & = \int_{\R}\sum_{n = 0}^{\infty}\dfrac{1+t^2}{(t-x_n)^2+y_n^2}d\mu_1(t)  \\
& \leq \int_{\R}\dfrac{2}{\tilde{\beta}}\dfrac{1+t^2}{\abs{t}+x_0-\tilde{\beta}/2}d\mu_1(t) \\
& \leq \int_{\R}C(1+\abs{t})d\mu(t) < + \infty
\end{align*}
for some $C > 0$, where we have used \eqref{eq:ineqsum} and the fact that the functions we are integrating are positive.
\end{proof}
\end{lemma}
We can now give the characterization for the two final cases:
\begin{proposition}
\label{prop:lastcase}
\begin{enumerate}[\hspace{0.5cm}\normalfont(a)]
\item Let $\tilde{\beta} > 0$ and suppose that $\mu$ is such that
$$\int_{\R}\abs{t}d\mu(t) < \infty, \quad \int_{(-\infty,0)}t^2d\mu(t) = \infty, \quad \int_{(0,+\infty)}t^2d\mu(t) < \infty.$$
Then $f$ is of PHS.
\item Let $\tilde{\beta} < 0$ and suppose that $\mu$ is such that
$$\int_{\R}\abs{t}d\mu(t) < \infty, \quad \int_{(-\infty,0)}t^2d\mu(t) < \infty, \quad \int_{(0,+\infty)}t^2d\mu(t) = \infty.$$
Then $f$ is of PHS.
\end{enumerate}
\begin{proof}
By Lemma \ref{lemma:eps1}, $y_n$ is bounded, and so Proposition \ref{prop:phsiffconv} assures that $f$ is of PHS.
\end{proof}
\end{proposition}

With the previous result, we have a complete characterization of the hyperbolic step of parabolic functions given by Theorem \ref{teo:expresion} when $t \in L^1(\mu)$. With the notation introduced in the beginning of Section \ref{sec:t1}, this characterization can be summarized in the following table

\begin{table}[H]
\renewcommand{\arraystretch}{1.3}
\begin{tabular}{|c|cc|c|c|}
\hline
\multirow{5}{*}{$\tilde{\beta}$} & \multicolumn{2}{c|}{\multirow{3}{*}{$\displaystyle\int_I t^2 d\mu(t)$}} & \multirow{5}{*}{Hyp. Step} & \multirow{5}{*}{Reference} \\
 & \multicolumn{2}{c|}{} & & \\
 & \multicolumn{2}{c|}{} & & \\ \cline{2-3}
 & \multicolumn{1}{c|}{\multirow{2}{*}{$I = (-\infty,0)$}} & \multirow{2}{*}{$I = (0,+\infty)$} & & \\ 
 & \multicolumn{1}{c|}{} & & & \\ \hline
$\tilde{\beta} = 0$    & \multicolumn{1}{c|}{$\ast$}    & $\ast$    & 0HS & Proposition \ref{prop:beta0}     \\ \hline
$\tilde{\beta} \neq 0$ & \multicolumn{1}{c|}{Converges} & Converges & PHS & Theorem \ref{teo:t2case}         \\ \hline
$\tilde{\beta} > 0$    & \multicolumn{1}{c|}{$\ast$}    & Diverges  & 0HS & Proposition \ref{prop:necessary} \\ \hline
$\tilde{\beta} > 0$    & \multicolumn{1}{c|}{Diverges}  & Converges & PHS & Proposition \ref{prop:lastcase}  \\ \hline
$\tilde{\beta} < 0$    & \multicolumn{1}{c|}{Converges} & Diverges  & PHS & Proposition \ref{prop:lastcase}  \\ \hline
$\tilde{\beta} < 0$    & \multicolumn{1}{c|}{Diverges}  & $\ast$    & 0HS & Proposition \ref{prop:necessary} \\ \hline
\end{tabular}
\end{table}

The content of this section can be related to \cite[Section 2]{CDP}, where the authors define the following:
\begin{definition}
A parabolic function $g \colon \D \to \D$ with DW point $\tau$ is said to be of angular-class of order $p \in \N$ at $\tau$, denoted as $g \in C^p_A(\tau)$, if
$$g(z) = \tau + \sum_{j=1}^p\dfrac{a_j}{j!}(z-\tau)^j+\gamma(z), \quad z \in \D,$$
where $a_1,\ldots,a_p \in \C$ and $\gamma$ is a holomorphic function on $\D$ with
$$\angle\lim_{z \to \tau}\dfrac{\gamma(z)}{(z-\tau)^p}=0.$$
\end{definition}
They also prove the following characterization:
\begin{proposition}
\cite[Propositions 2.1 and 2.2]{CDP}
Let $g \colon \D \to \D$ be a parabolic function whose DW point is $\tau$ and consider its conjugation $f \colon \H \to \H$ whose DW point is infinity.
\begin{enumerate}[\hspace{0.5cm}\normalfont(a)]
\item $g \in C^2_A(\tau)$ if and only if $\angle\lim_{z \to \infty}(f(z)-z) \in \R \cup \H$.
\item If $g \in C^2_A(\tau)$, then $g \in C^3_A(\tau)$ if and only if there exists $a,b \in \C$ with $\angle\lim_{z \to \infty}z(f(z) - z - a - b/z) = 0$.
\end{enumerate}
\end{proposition}
From this, one can see that Lemma \ref{lemma:anglep0} implies that functions satisfying the condition $t \in L^1(\mu)$ can be conjugated to a function in $C^2_A(\tau)$, for every $\tau \in \partial\D$. Similarly, Proposition \ref{prop:t2charac} implies that functions satisfying the condition $t \in L^2(\mu)$ can be conjugated to a function in $C^3_A(\tau)$. However, the map $z \in \H \mapsto z + i \in \H$ is related to both classes $C^2_A(\tau)$ and $C^3_A(\tau)$, but it does not satisfy the condition that $t \in L^1(\mu)$ (it comes from $\beta = 0$ and $d\mu/dm(t) = 1/(1+t^2)$, where $m$ is the Lebesgue's measure on $\R$).

\section{Characterization when $\mu$ is symmetric}
\label{sec:symmetric}
In this section, instead of considering integrability conditions on $\mu$, we focus on symmetry. The measure $\mu$ in Theorem \ref{teo:expresion} is said to be symmetric if $\mu(A) = \mu(-A)$ for every measurable set $A \subset \R$. This property was already considered by Aaronson in \cite[Page 217]{Aaronson}. In this section, we will suppose that $f$ is given by the expression in Theorem \ref{teo:expresion}, and the corresponding measure $\mu$ is symmetric. In this setting, we will try to characterize the hyperbolic step of $f$ in terms of $\beta$.

For example, if the associated measure of $f$ is symmetric and $\beta = 0$, then it is possible to check that $f(\{it : t > 0\}) \subset \{it : t > 0\}$. Then, every orbit of $f$ converges non-tangentially to infinity, and so the following known result is derived (see Theorem \ref{teo:convtang}):
\begin{proposition}
\label{prop:AaronsonSymmetric}
If $\mu$ is symmetric and $\beta = 0$, then $f$ is of 0HS.
\end{proposition}

As before, the case of $\beta \neq 0$ is more subtle. To tackle this problem, let us rewrite $f$ as
$$f(z) = z + \beta + p_0(z) - \dfrac{\mu(\{0\})}{z}, \quad z \in \H,$$
where
$$p_0(z) = 2z\int_{(0,+\infty)}\dfrac{1+t^2}{t^2-z^2}d\mu(t).$$
Indeed, the orbits satisfy
$$\left\lbrace\begin{array}{l}
x_{n+1} = x_n + \beta + x_np_1(z_n)-\text{Re}(\mu(\{0\})/z_n), \\
y_{n+1} = y_n(1+p_2(z_n))-\text{Im}(\mu(\{0\})/z_n),
\end{array}
\right.$$
where
$$p_1(z) = \int_{(0,+\infty)}\dfrac{2(1+t^2)(t^2-\abs{z}^2)}{\abs{t^2-z^2}^2}d\mu(t),$$
$$p_2(z) = \int_{(0,+\infty)}\dfrac{2(1+t^2)(t^2+\abs{z}^2)}{\abs{t^2-z^2}^2}d\mu(t).$$
In particular,
$$\abs{p_1(z)} \leq p_2(z), \quad z \in \H.$$
As said before, we will suppose that
$$\int_{\R}\abs{t}d\mu(t) = \infty,$$
as in other case the function falls in the hypothesis of Theorem \ref{teo:t1end}. With this setting, the following is satisfied:
\begin{lemma}
\label{lemma:yndivergence}
If $\mu$ is symmetric with $t \not\in L^1(\mu)$ and $\beta \neq 0$, then $\lim_{n \to \infty} y_n = +\infty$ for any $z_0 \in \H$.
\begin{proof}
For brevity, denote $p_2(z_n) = \varepsilon_n$. Observe that
$$\abs{z_{n+1}} \leq \abs{z_n}+\abs{\beta}+\abs{p_0(z_n)}+\dfrac{\mu(\{0\})}{\abs{z_n}} \leq \abs{\beta} + (1+\varepsilon_n)\abs{z_n} +\dfrac{\mu(\{0\})}{\abs{z_n}}.$$

Also, $\abs{z_n} \geq y_n \geq y_0$, and so
$$\abs{z_{n+1}} \leq C+\abs{\beta}+(1+\varepsilon_n)\abs{z_n},$$
where
$$C = \dfrac{\mu(\{0\})}{y_0}.$$

As $y_n$ is a non-decreasing sequence, if $\lim_n y_n \neq + \infty$, then
$$\lim_{n \to \infty}y_n = L \in (0,+\infty).$$
Note that
$$\dfrac{y_{n+1}}{y_n} = 1+\varepsilon_n+\dfrac{\mu(\{0\})}{\abs{z_n}^2} \geq 1+\varepsilon_n.$$
Then,
$$y_{n+1} \geq y_0\prod_{j = 1}^n(1+\varepsilon_j).$$
In that case, assuming $y_n$ converges, we have
$$\prod_{j = 1}^{\infty}(1+\varepsilon_j) < \infty.$$

Using the former ideas, one has
\begin{align}
\label{eq:bound1}
\abs{z_{n+1}} & \leq C+\abs{\beta}+(1+\varepsilon_n)\abs{z_n} \leq (1+\varepsilon_n)(C+\abs{\beta}+\abs{z_{n}}) \\
& \leq (1+\varepsilon_n)(2(C+\abs{\beta})+(1+\varepsilon_{n-1})\abs{z_{n-1}}) \nonumber \\
& \leq (1+\varepsilon_n)(1+\varepsilon_{n-1})(2(C+\abs{\beta})+\abs{z_{n-1}}) \leq \cdots \nonumber \\
& \leq \prod_{j = m}^n(1+\varepsilon_j)((n-m)(C+\abs{\beta})+\abs{z_m}), \quad n > m \geq 0. \nonumber
\end{align}

Indeed, by convergence, it is possible to find $n_0 \in \N$ such that
\begin{equation}
\label{eq:bound2}
\prod_{j = n_0}^{\infty}(1+\varepsilon_j)  < 2.    
\end{equation}

Then, for $n > m \geq n_0$, we have
$$\abs{z_n} \leq 2(n-m)(C+\abs{\beta})+2\abs{z_m}.$$

Also, denote
$$I_n = \left(\dfrac{1}{2}\abs{z_n},2\abs{z_n}\right).$$
Note that there exists $K > 0$ such that if $t \in I_n$, then
$$\dfrac{2(1+t^2)(t^2+\abs{z_n}^2)}{\abs{t^2-z_n^2}^2} \geq K.$$
In that case, $\varepsilon_n \geq K\mu(I_n)$. On the other hand, as
$$\prod_{j = 1}^{\infty}(1+\varepsilon_j) < \infty,$$
then
$$\sum_{n = 1}^{\infty}\varepsilon_n < \infty.$$
So, we arrive to a contradiction if
$$\sum_{n = 1}^{\infty} \mu(I_n) = \infty.$$

Let $g$ be the function given by
$$g(t) = \sum_{n = 1}^{\infty}\chi_{I_n}(t) = \#\{n \in \N \mid t \in I_n\}, \quad t > 0.$$
Then,
$$\sum_{n = 1}^{\infty} \mu(I_n) = \int_{(0,+\infty)}g(t)d\mu(t).$$
We claim that
$$\int_{(0,+\infty)}g(t)d\mu(t) = \infty,$$
getting the desired contradiction.

To prove the claim, let $t_0 > 0$ be such that $\max\{\abs{z_0},\ldots,\abs{z_{n_0}}\} < t_0/2$, where $n_0$ is defined by means of \eqref{eq:bound2}. Given $t > t_0$, define $m(t) = \max\{k \in \N \mid \abs{z_k} < t/2\}$. Note that if $t > t_0$, then $m(t) \geq n_0$. Also, if $n > m(t)$, then $\abs{z_n} > t/2$. In that case, let $n > m(t)$ be such that
$$2(C+\abs{\beta})(n-m(t)) \leq t.$$
With this, we have
$$\abs{z_n} \leq 2(n-m(t))(C+\abs{\beta})+2\abs{z_{m(t)}} \leq t+t=2t.$$
These calculations show that for $t > t_0$ we have $t \in I_n$ at least when
$$n \geq m(t), \quad 2(C+\abs{\beta})(n-m(t)) \leq t.$$

Joining all these ideas, if $t > t_0$, then
\begin{align*}
g(t) & = \#\{n \in \N \mid t \in I_n\} \\
& \geq \#\{n \in \N \mid n > m(t), \, 2(C+\abs{\beta})(n-m(t)) \leq t\} \geq \dfrac{1}{2(C+\abs{\beta})}t-1.
\end{align*}
So, 
$$\int_{(0,+\infty)}g(t)d\mu(t) \geq \dfrac{1}{2(C+\abs{\beta})}\int_{(t_0,+\infty)}td\mu(t) - \mu(\R) = +\infty.$$
\end{proof}
\end{lemma}

We give now some properties of the orbits:
\begin{lemma}
\label{lemma:xCy}
If $\mu$ is symmetric with $t \not\in L^1(\mu)$ and $\beta \neq 0$, then 
$$\abs{x_n} \leq \left(\abs{\beta}+\dfrac{\mu(\{0\})}{y_0}\right)y_n\left(\dfrac{1}{y_n} + \cdots + \dfrac{1}{y_1}+\dfrac{\abs{x_0}}{y_0}\right).$$
\begin{proof}
For brevity, denote $\varepsilon_n = p_2(z_n)$. Note that
$$\abs{x_n} \leq \abs{\beta}+(1+\varepsilon_{n-1})\abs{x_{n-1}}+\dfrac{\mu(\{0\})}{\abs{y_0}} = \left(\abs{\beta}+\dfrac{\mu(\{0\})}{y_0}\right)+(1+\varepsilon_{n-1})\abs{x_{n-1}}.$$
But
$$y_n = y_{n-1}(1+\varepsilon_{n-1})-\text{Im}(\mu(\{0\})/z_n) \geq y_{n-1}(1+\varepsilon_{n-1}),$$
because $-\text{Im}(\mu(\{0\})/z_n) > 0$. In that case
$$\abs{x_n} \leq \left(\abs{\beta}+\dfrac{\mu(\{0\})}{y_0}\right)+\dfrac{y_n}{y_{n-1}}\abs{x_{n-1}}.$$

Denote
$$C = \abs{\beta}+\dfrac{\mu(\{0\})}{y_0}, \quad \alpha_n = \dfrac{y_n}{y_{n-1}}.$$
Then, we just showed $\abs{x_n} \leq C + \alpha_n\abs{x_{n-1}}$. Applying this inequality repeatedly, one has
\begin{align*}
\abs{x_n} & \leq C + \alpha_n\abs{x_{n-1}} \leq C + \alpha_n(C+\alpha_{n-1}\abs{x_{n-2}}) \\
& = C(1+\alpha_n) + \alpha_n\alpha_{n-1}\abs{x_{n-2}} \\
& \leq C(1+\alpha_n+\alpha_n\alpha_{n-1})+\alpha_n\alpha_{n-1}\alpha_{n-2}\abs{x_{n-3}} \\
& \leq C(1+\alpha_n+\alpha_n\alpha_{n-1}+\cdots+\alpha_n\cdots\alpha_2)+\alpha_n\cdots\alpha_1\abs{x_0} \\
& = C\left(1+\sum_{j = 2}^n\prod_{k=j}^n\alpha_k\right)+\prod_{k=1}^n\alpha_k\abs{x_0}.
\end{align*}
Note that, by definition,
$$\prod_{k=j}^n\alpha_k = \dfrac{y_n}{y_{j-1}}.$$
Then, we can rewrite
\begin{align*}
\abs{x_n} & \leq  C\left(1+y_n\sum_{j = 1}^{n-1}\dfrac{1}{y_j}\right)+\dfrac{y_n}{y_0}\abs{x_0} = Cy_n\left(\dfrac{1}{y_n} + \cdots + \dfrac{1}{y_1}+\dfrac{\abs{x_0}}{y_0}\right).
\end{align*}
\end{proof}
\end{lemma}
\begin{lemma}
\label{lemma:on}
If $\mu$ is symmetric with $t \not\in L^1(\mu)$ and $\beta \neq 0$, then $\displaystyle\lim_{n \to \infty}\dfrac{\abs{x_n}}{ny_n} = 0$.
\begin{proof}
It is a mere consequence via Stolz's criterion. Denote
$$A_n = \left(\abs{\beta}+\dfrac{\mu(\{0\})}{y_0}\right)\left(\dfrac{1}{y_n} + \cdots + \dfrac{1}{y_1}+\dfrac{\abs{x_0}}{y_0}\right).$$
Then, by Lemma \ref{lemma:xCy},
\begin{align*}
\lim_{n \to \infty}\dfrac{\abs{x_n}}{ny_n} & \leq \lim_{n \to \infty}\dfrac{A_n}{n} = \lim_{n \to \infty}\dfrac{A_{n+1}-A_n}{(n+1)-n} \\
& = \lim_{n \to \infty}\left(\abs{\beta}+\dfrac{\mu(\{0\})}{y_0}\right)\dfrac{1}{y_{n+1}} = 0,
\end{align*}
due to Lemma \ref{lemma:yndivergence}.
\end{proof}
\end{lemma}
\begin{lemma}
\label{lemma:liminf}
If $\mu$ is symmetric with $t \not\in L^1(\mu)$ and $\beta \neq 0$, then 
$$\liminf_{n \to \infty}\dfrac{x_n}{y_n}p_2(z_n) = 0.$$
\begin{proof}
Suppose that there exist $\delta > 0$ and $N_1 \in \N$ such that
$$\dfrac{\abs{x_n}}{y_n}p_2(z_n) \geq \delta, \quad n \geq N_1.$$
Applying Lemma \ref{lemma:on}, choose $N_2 \in \N$ such that
$$\dfrac{\abs{x_n}}{ny_n} \leq \dfrac{\delta}{4}, \quad n \geq N_2.$$
Denote $N = \max\{N_1,N_2\}$. Then, $p_2(z_n) \geq 4/n$ for every $n \geq N$. Also,
$$y_n = (1+p_2(z_{n-1}))y_{n-1}-\text{Im}\left(\dfrac{\mu(\{0\})}{z_{n-1}}\right) \geq \left(1+\dfrac{4}{n-1}\right)y_{n-1}, \quad n > N.$$
Applying this argument several times, one has
$$y_n \geq \prod_{k = N}^{n-1}\left(1+\dfrac{4}{k}\right)y_N, \quad n > N.$$

Observe that
$$1+\dfrac{1}{x} \geq \exp\left(\dfrac{1}{2x}\right), \quad x \geq 1.$$
In that case,
$$y_n \geq \exp\left(\sum_{k = N}^{n-1}\dfrac{2}{k}\right)y_N, \quad n \geq N,$$
and so
$$\lim_{n \to \infty} \dfrac{y_n}{n^2} \geq \lim_{n \to \infty} \exp\left(\sum_{k = N}^{n-1}\dfrac{2}{k}\right)\dfrac{y_N}{n^2} = Ky_N,$$
for some $K = K(N) > 0$. Therefore, it is clear that
$$\sum_{n = 0}^{\infty}\dfrac{1}{y_n} < \infty.$$
Then, by Lemma \ref{lemma:xCy}, it is possible to find $C > 0$ such that $\abs{x_n} \leq Cy_n$. On the other hand,
$$p_2(z_n) = \dfrac{y_{n+1}-y_n}{y_n}+\dfrac{1}{y_n}\text{Im}\left(\dfrac{\mu(\{0\})}{z_n}\right) \to 0, \quad n \to \infty,$$
as both terms converge to zero (for the first of them, see Theorem \ref{teo:Pomm}). Finally,
$$\liminf_{n \to \infty}\dfrac{\abs{x_n}}{y_n}p_2(z_n) \leq C\liminf_{n \to \infty}p_2(z_n) = 0,$$
This leads to a contradiction.
\end{proof}
\end{lemma}

We can now prove the following:
\begin{theorem}
\label{teo:symmetric}
If $\mu$ is symmetric with $t \not\in L^1(\mu)$, then $f$ is of 0HS for every $\beta \in \R$.
\begin{proof}
If $\beta = 0$, apply Proposition \ref{prop:AaronsonSymmetric}. In other case, note that
$$\dfrac{x_{n+1} - x_n}{y_n} = \dfrac{\beta + p_1(z_n)x_n-\text{Re}(\mu(\{0\})/z_n)}{y_n}.$$
First, $\text{Re}(\mu(\{0\})/z_n) \to 0$ as $n \to \infty$, because $z_n \to \infty$. On the other hand, applying Lemma \ref{lemma:yndivergence}, one has
$\beta/y_n \to 0$ as $n \to \infty$. So, by Lemma \ref{lemma:liminf},
\begin{align*}
\lim_{n \to \infty} \dfrac{\abs{x_{n+1} - x_n}}{y_n} & = \liminf_{n \to \infty}  \dfrac{\abs{x_{n+1} - x_n}}{y_n} = \liminf_{n \to \infty} \dfrac{\abs{x_n}}{y_n}\abs{p_1(z_n)} \\
& \leq \liminf_{n \to \infty} \dfrac{\abs{x_n}}{y_n}p_2(z_n) = 0. 
\end{align*}
Then, $f$ is of 0HS, as Theorem \ref{teo:Pomm} assures.
\end{proof}
\end{theorem}

With the previous result, we have a complete characterization of the hyperbolic step of parabolic functions given by Theorem \ref{teo:expresion} when $\mu$ is a symmetric measure. With the usual notation, this characterization can be summarized in the following table

\begin{table}[H]
\centering
\renewcommand{\arraystretch}{1.3}
\begin{tabular}{|c|c|c|c|}
\hline
$\beta$        & $\mu$                               & Hyp. Step & Reference \\ \hline
$\beta = 0$    & *                                   & 0HS       & Proposition \ref{prop:AaronsonSymmetric} \\ \hline
$\beta \neq 0$ & $t \in L^2(\mu)$                    & PHS       & Theorem \ref{teo:t2case} \\ \hline 
$\beta \neq 0$ & $t \in L^1(\mu) \setminus L^2(\mu)$ & 0HS       & Proposition \ref{prop:necessary} \\ \hline 
$\beta \neq 0$ & $t \not\in L^1(\mu)$                & 0HS       & Theorem \ref{teo:symmetric} \\ \hline
\end{tabular}

\end{table}
\section{Characterization when $\mu$ is (mainly) supported on a half-line}
\label{sec:perturb}
This section is devoted to develop a proof of Theorem \ref{teo:perturb}, which will be shown in Subsection \ref{subsec:main}. This proof uses a construction which is based on the behaviour of the orbits of parabolic functions whose related measure satisfies more restrictive hypothesis than the ones in Theorem \ref{teo:perturb}. To introduce such topic, let $f$ be given by the expression of Theorem \ref{teo:expresion}, and suppose that $\mu$ is supported on a half-line. That is, assume that there exists $a \in \R$ with $\mu((a,+\infty)) = 0$.
As in the previous section, having in mind the results of Section \ref{sec:t1}, we will focus on the case
$$\int_{\R}\abs{t}d\mu(t) = \infty.$$

Our main interest is to characterize the hyperbolic step of $f$ under such conditions. In order to do so, let us define $p \colon \C \setminus (-\infty,a] \to \C$ by
$$p(z) = \int_{(-\infty,a]}\dfrac{1+tz}{t-z}d\mu(t), \quad z \in \C\setminus(-\infty,a].$$
Note that $p$ is holomorphic, and then $f$ can be extended to a holomorphic map $\hat{f} \colon \C\setminus(-\infty,a] \to \C$ with
$$\hat{f}(z) = z + \beta + p(z), \quad z \in \C\setminus(-\infty,a].$$
The next lemma summarizes some properties of $\hat{f}$.
\begin{lemma}
\label{lemma:properties}
\begin{enumerate}[\normalfont(a)]
\item $\hat{f}(\H) \subset \H$, $\hat{f}((a,+\infty)) \subset \R$, $\hat{f}(-\H) \subset -\H$.
\item $\textnormal{Re}(p(x))$ is an increasing function for $x > a$.
\item $\textnormal{Re}(p(x+iy)) \geq \textnormal{Re}(p(x))$ for any $x > a$, $y \in \R$.
\item $\displaystyle\lim_{x \to + \infty}\textnormal{Re}(p(x)) = + \infty$.
\item There exists $b > a$ such that $p(b) + \beta > 0$. Consequently, $\hat{f}(\C_b) \subset \C_b$, where $\C_b = \{z \in \C : \textnormal{Re}(z) > b\}$.
\item For every $b > a$ such that $\hat{f}(\C_b) \subset \C_b$, both the map $\hat{f} \colon \C \setminus (-\infty,b] \to  \C \setminus (-\infty,b]$ and the restricted map $g = \hat{f}|_{\C_b} \colon \C_b \to \C_b$ are parabolic of 0HS, and their DW point is infinity.
\end{enumerate}
\begin{proof}
(a) Note that
$$\text{Im}(\hat{f}(x+iy)) = y\left(1+\int_{\R}\dfrac{1+t^2}{(t-x)^2+y^2}d\mu(t)\right),$$
and so $\text{Im}(\hat{f}(x+iy))$ has the same sign as $y$.

(b) Since
$$\text{Re}(p(x)) = \int_{(-\infty,a]}\dfrac{1+tx}{t-x}d\mu(t),$$
and
$$\dfrac{\partial}{\partial x}\left(\dfrac{1+tx}{t-x}\right) = \dfrac{1+t^2}{(t-x)^2} > 0, \quad t \leq a < x,$$
the statement holds.

(c) Note that
$$\text{Re}(p(x+iy)) = \int_{\R}\dfrac{t-x+t^2x-t(x^2+y^2)}{(t-x)^2+y^2}d\mu(t),$$
and
$$\dfrac{\partial}{\partial y}\left(\dfrac{t-x+t^2x-t(x^2+y^2)}{(t-x)^2+y^2}\right) = -\dfrac{2y(1+t^2)(t-x)}{((t-x)^2+y^2)^2} \geq 0,$$
whenever $t \leq a < x$, $y \geq 0$. Thus, (c) is clear.

(d) We need to show that
$$\lim_{x \to +\infty}\int_{(-\infty,a]}\dfrac{1+tx}{t-x}d\mu(t) = + \infty.$$
To do so, note that
$$\dfrac{1+tx}{t-x} > 0 \iff t < -\dfrac{1}{x},$$
where we are considering that $t \leq a < x$ and $x > 0$. By Fatou's lemma, one has
$$\lim_{x \to +\infty}\int_{(-\infty,-1/x)}\dfrac{1+tx}{t-x}d\mu(t) \geq \int_{(-\infty,0)}(-t)d\mu(t) = +\infty.$$
If $a < 0$, we can consider $x$ large enough so that $(-\infty,a] \subset (-\infty, -1/x)$, and it is enough to apply the previous argument to conclude the proof. On the other hand, if $a \geq 0$, it is possible to check that
\begin{align*}
\int_{(-1/x,a)}\abs{\dfrac{1+tx}{t-x}}d\mu(t) & \leq \abs{\dfrac{1+ax}{a-x}}\mu((-1/x,a)) \\
& \to a\mu([0,a)) < \infty, \quad x \to + \infty.
\end{align*}

Joining both arguments, whether $a < 0$ or $a \geq 0$, one has
$$\lim_{x \to + \infty} \int_{\R}\dfrac{1+tx}{t-x}d\mu(t) = +\infty.$$

(e) Using previous properties, it is enough to find $b > a$ such that
$$\hat{f}(x)-x = p(x) +\beta > 0, \quad x > b,$$
which can be done by (d).

(f) Let us conjugate $\hat{f}$ to the function $F \colon \C_0 \to \C_0$ given by
$$F(z) = \sqrt{\hat{f}(z^2+b)-b}, \quad z \in \C_0,$$
where we are using the principal argument to define the square root. One can see that
\begin{align*}
\angle\lim_{z \to \infty}\dfrac{F(z)}{z} & = \angle\lim_{z \to \infty}\sqrt{\dfrac{\hat{f}(z^2+b)-b}{z^2}} \\
& = \angle\lim_{z \to \infty}\sqrt{\dfrac{\hat{f}(z^2+b)-b}{z^2+b}\dfrac{z^2+b}{z^2}} = 1,
\end{align*}
because
$$\angle\lim_{z \to \infty}\dfrac{f(z)}{z} = 1.$$
This shows that $F$ is parabolic, and so is $\hat{f}$. Also, both $F$ and $\hat{f}$ have infinity as their DW point. Moreover, it is of 0HS, because $\hat{f}((b,+\infty)) \subset (b,+\infty)$, and so any orbit whose initial point lies in this interval converges non-tangentially (see Theorem \ref{teo:convtang}). Similar ideas can be used in the case of $g$, under the conjugation $z \in \C_b \mapsto z-b \in \C_0$.
\end{proof}
\end{lemma}

The proof of the Theorem \ref{teo:perturb} can be achieved by joining several related ideas. In particular, in order to be as clear as possible, we will start by showing a weaker version of such theorem:
\begin{proposition}
\label{prop:weak}
Suppose that there exists $a \in \R$ such that $\mu((a,+\infty)) = 0$ and
$$\int_{\R}\abs{t}d\mu(t) = +\infty.$$
Then, $f$ is of PHS for every $\beta \in \R$.
\end{proposition}
Indeed, two different proofs of this result will be shown. In Subsection \ref{subsec:firstproof} we give a proof based on the existence of the Koenigs function built in \cite{BakerPommerenke}. In Subsection \ref{subsec:secondproof} we give a proof based on a characterization of the hyperbolic step by Pommerenke, see Theorem \ref{teo:Pomm}. Later, in Subsection \ref{subsec:main}, we prove Theorem \ref{teo:perturb} by adapting the ideas in Subsection \ref{subsec:secondproof}.

\subsection{First proof of Proposition \ref{prop:weak}}
\label{subsec:firstproof}
Let $b > a$ and $g$ be as in Lemma \ref{lemma:properties}, define $z_0 = b+1$, and consider $z_n = g^n(z_0)$ for $n \in \N$. In \cite{BakerPommerenke}, it is shown that the functions $h_n \colon \C_b \to \C$ given by
$$h_n(z) = \dfrac{g^n(z)-z_n}{z_{n+1}-z_n}, \quad z \in \C_b,$$
converge uniformly on compact sets to a function $h \colon \C_b \to \C$ satisfying
$$h \circ g = h + 1.$$
Note that it has been used that the translation $z \in \C_0 \mapsto z+b \in \C_b$ is a bijective map, and so $g$ can be conjugated to a function from $\C_0$ to $\C_0$, as usual.

Since $g((b,+\infty)) \subset (b,+\infty)$, then $b < z_n < z_{n+1}$. In particular, for every $n \in \N$ one has
$$h_n(\C_b \cap \H) \subset \H, \quad h_n((b,+\infty)) \subset \R, \quad h_n(\C_b \cap (-\H)) \subset \C_b \cap (-\H).$$
Taking the limit as $n \to +\infty$, because $h$ is an open map (note that it can not be constant because $h \circ g = h + 1$), it also satisfies
$$h(\C_b \cap \H) \subset \H, \quad h((b,+\infty)) \subset \R, \quad h(\C_b \cap (-\H)) \subset \C_b \cap (-\H).$$

The function $h$ can be extended to a function $\hat{h} \colon \C \setminus (-\infty,b] \to \C$ in the following way: let $x > b$, and denote the pseudo-hyperbolic distance in $\C \setminus (-\infty,b]$ by $\tilde{\rho}$. Note that
\begin{align*}
\min_{y \neq 0}\tilde{\rho}(x,b+iy) & = \min_{y \neq 0}\rho_{\C_0}(\sqrt{x-b},\sqrt{iy}) = \min_{r > 0}\rho_{\C_0}(1,re^{\pm i\pi/4}) \\
& = \rho_{\C_0}(1,e^{\pm i\pi/4}) = \sqrt{2}-1,
\end{align*}
where
$$\rho_{\C_0}(z,w) = \abs{\dfrac{z-w}{z+\overline{w}}}, \quad z,w \in \C_0.$$
With that in mind, let $w_0 \in \C \setminus (-\infty,b]$. By Lemma \ref{lemma:properties} and \cite[Corollary 4.6.9.(iv)]{AbateNewBook}, $\tilde{\rho}(\hat{f}^n(z_0),\hat{f}^n(w_0)) \to 0$, and so it is possible to find $N = N(w_0) \in \N$ such that
$$\tilde{\rho}(\hat{f}^n(z_0),\hat{f}^n(w_0)) < \sqrt{3-2\sqrt{2}}, \quad n \geq N.$$
Then, $f^n(w_0) \in \C_b$ for every $n \in \N$ with $n \geq N$. Similarly, if $K \subset \C \setminus (-\infty,b]$ is a compact set, it is also possible to find $N = N(K) \in \N$ such that $f^n(w) \in \C_b$ for every $n \in \N$ with $n \geq N$ and for every $w \in K$. Indeed, if $n,m \in \N$ are such that $n > m \geq N(K)$,  
$$h(\hat{f}^n(z)) = h(\hat{f}^{n-m}(\hat{f}^m(z))) = h(g^{n-m}(\hat{f}^m(z))) = h(\hat{f}^m(z)) + n - m,$$
and so
$$h(\hat{f}^n(z)) - n = h(\hat{f}^m(z)) - m.$$
Using this property, one can define the function $\hat{h} \colon \C \setminus (-\infty,b] \to \C$ by
$$\hat{h}(z) = h(\hat{f}^n(z))-n, \quad z \in \C \setminus (-\infty,b],$$
where $n \in \N$ is such that $\hat{f}^n(z) \in \C_b$, as above. Due to the previous discussion, $\hat{h}$ is well-defined, and it is a holomorphic function (it can be defined locally as a composition of holomorphic maps in every compact set).

Furthermore, if $z \in \H$ and $n \in \N$ is such that $\hat{f}^n(z) \in \C_b$, then $\hat{h}(z) = h(\hat{f}^n(z))-n \in \H$, because $h(\C_b \cap \H) \subset \H$. Thus, $\hat{h}(\H) \subset \H$. Moreover, by analytic continuation, 
$$\hat{h} \circ \hat{f} = \hat{h}+1, \quad z \in \C \setminus (-\infty,a].$$

Let us now define the restriction $\sigma = \hat{h}|_{\H} \colon \H \to \H$. Note that if $z \in \H$, then
$$\sigma(f(z)) = \hat{h}(\hat{f}(z)) = \hat{h}(z)+1 = \sigma(z)+1.$$
In particular, for any $z_0 \in \H$, Schwarz-Pick Lemma can be used to show that
\begin{align*}
\rho_{\H}(z_n,z_{n+1}) & \geq \rho_{\H}(\sigma(z_n),\sigma(z_{n+1})) = \rho_{\H}(\sigma(z_0)+n,\sigma(z_0)+n+1) \\
& = \rho_{\H}(\sigma(z_0),\sigma(z_0)+1) = C(z_0) > 0,
\end{align*}
and so $f$ is of PHS.

\subsection{Second proof of Proposition \ref{prop:weak}}
\label{subsec:secondproof}
Let $b > a$ as in Lemma \ref{lemma:properties}. For $x > b$, $y \geq 0$, define the functions given by
\begin{align*}
P(x,y) = \beta + \int_{(-\infty,a]}\dfrac{t-x+t^2x-t(x^2+y^2)}{(t-x)^2+y^2}d\mu(t)
\end{align*}
and
$$Q(x,y) = \displaystyle\int_{(-\infty,a]}\dfrac{1+t^2}{(t-x)^2+y^2}d\mu(t).$$
Note that, if $y > 0$, then
$$P(x,y) = \text{Re}(f(x+iy))-x, \quad Q(x,y) = \dfrac{\text{Im}(f(x+iy))-y}{y}.$$
Thus, if $z_0 \in \H$, $x_{n+1} = x_n + P_n$, where $P_n = P(x_n,y_n)$, and $y_{n+1} = y_n(1+Q_n)$, where $Q_n = Q(x_n,y_n)$.

Similarly, given $\tilde{x}_0 > b$, define $\tilde{x}_{n+1} = \hat{f}(\tilde{x}_n) = \tilde{x}_n + \tilde{P}_n$, where $\tilde{P}_n = P(\tilde{x}_n,0)$. Also, given $\tilde{y}_0 > 0$, define $\tilde{y}_{n+1} = \tilde{y}_n(1+\tilde{Q}_n)$, where $\tilde{Q}_n = Q(\tilde{x}_n,0)$. The quantities $\tilde{x}_n$ and $\tilde{y}_n$ might be understood as the boundary counterparts of $x_n$ and $y_n$. Indeed, the following relations among them are satisfied:
\begin{lemma}
\label{lemma:secondproof}
Suppose that $b < \tilde{x}_0 \leq x_0$ and $\tilde{y}_0 \geq y_0 > 0$. Then, for all $n \in \N$,
\begin{enumerate}[\hspace{0.5cm}\normalfont(a)]
\item $0 < \tilde{P}_n \leq P_n$ and $\tilde{x}_n \leq x_n$,
\item $\tilde{y}_n \geq y_n$,
\item $\tilde{y}_{n+1} \geq \tilde{y}_n$ and $\tilde{P}_n \leq \tilde{P}_{n+1}$.
\end{enumerate}
\begin{proof}
(a) Observe that
$$\tilde{P}_n  = \beta + p(\tilde{x}_n) \geq \beta + p(b) > 0.$$

Also, note that
$$\dfrac{\partial P}{\partial y}(x,y) = \int_{(-\infty,a]}-\dfrac{2y(1+t^2)(t-x)}{((t-x)^2+y^2)^2}d\mu(t).$$
In particular, if $x > a$, then
$$\dfrac{\partial P}{\partial y}(x,y) > 0.$$

On the other hand,
$$\dfrac{\partial P}{\partial x}(x,0) = \int_{(-\infty,a]}\dfrac{1+t^2}{(t-x)^2}d\mu(t).$$
So, if $x > a$, then
$$\dfrac{\partial P}{\partial x}(x,0) > 0.$$

Using these properties, the proof can be done by induction. Suppose that $\tilde{x}_n \leq x_n$ and $\tilde{P}_n \leq P_n$. In that case,
$$\tilde{x}_{n+1} = \tilde{x}_n+\tilde{P}_n \leq x_n + P_n = x_{n+1}.$$
Furthermore, using the monotonicity,
$$\tilde{P}_{n+1} = P(\tilde{x}_{n+1},0) \leq P(x_{n+1},0) \leq P(x_{n+1},y_{n+1}) = P_{n+1}.$$

(b) It is easy to see that $Q(x,y)$ is a decreasing function with respect to $x > a$. Similarly, it is also decreasing with respect to $y > 0$. So,
$$Q_n = Q(x_n,y_n) \leq Q(x_n,0) \leq Q(\tilde{x}_n,0) = \tilde{Q}_n.$$
Again, by induction, if $y_n \leq \tilde{y}_n$, then
$$y_{n+1} = y_n(1+Q_n) \leq \tilde{y}_n(1+\tilde{Q}_n) = \tilde{y}_{n+1}.$$

(c) Note that $\tilde{Q}_n \geq 0$, so $\tilde{y}_{n+1} \geq \tilde{y}_n$. Also,
$$\tilde{P}_n = P(\tilde{x}_n,0) \leq P(\tilde{x}_{n+1},0) = \tilde{P}_{n+1},$$
due to the monotonicity properties of $P$.
\end{proof}
\end{lemma}

With these ideas, in order to prove Proposition \ref{prop:weak}, pick $z_0 \in \H$ with $\text{Re}(z_0) > b$, as defined in Lemma \ref{lemma:properties}. The limit
$$L = \lim_{n \to + \infty}\dfrac{x_{n+1}-x_n}{y_n}$$
exists, as seen in Theorem \ref{teo:Pomm}. But, as seen in Lemma \ref{lemma:secondproof}, $x_{n+1}-x_n > 0$, and so $L \geq 0$. Remember that $f$ is of PHS if and only if $L > 0$, that is, if and only if there exists some constant $C > 0$ such that
$$\dfrac{y_n}{x_{n+1}-x_n} \leq C.$$
But,
$$\dfrac{y_n}{x_{n+1}-x_n} = \dfrac{y_n}{P_n} \leq \dfrac{\tilde{y}_n}{\tilde{P}_n}.$$
Thus, define
$$\tilde{C}_n = \dfrac{\tilde{y}_n}{\tilde{P}_n},$$
and note that $\tilde{C}_{n+1}/\tilde{C}_n \geq 1$ (by Lemma \ref{lemma:secondproof}.(c)) and
\begin{align*}
\dfrac{\tilde{C}_{n+1}}{\tilde{C}_n} - 1 & = \dfrac{\tilde{y}_{n+1}}{\tilde{y}_n}\dfrac{\tilde{P}_n}{\tilde{P}_{n+1}}-1 = (1+\tilde{Q}_n)\dfrac{\tilde{P}_n}{\tilde{P}_{n+1}}-1 = \dfrac{\tilde{Q}_n-(\tilde{P}_{n+1}-\tilde{P}_n)/\tilde{P}_n}{1+(\tilde{P}_{n+1}-\tilde{P}_n)/\tilde{P}_n}.
\end{align*}

But
$$\dfrac{\partial P}{\partial x}(x,0) = Q(x,0),$$
and $Q(x,0)$ is a decreasing function of $x$ (provided that $x > a$), so
$$\tilde{P}_{n+1}-\tilde{P}_n = \int_{\tilde{x}_n}^{\tilde{x}_{n+1}}Q(t,0)dt \geq (\tilde{x}_{n+1}-\tilde{x}_n)Q(\tilde{x}_{n+1},0) = \tilde{P}_n\tilde{Q}_{n+1}.$$
Therefore,
$$\dfrac{\tilde{C}_{n+1}}{\tilde{C}_n} - 1 \leq \dfrac{\tilde{Q}_n-\tilde{Q}_{n+1}}{1+\tilde{Q}_{n+1}} \leq \tilde{Q}_n-\tilde{Q}_{n+1}.$$
In particular,
$$\sum_{n = 0}^{\infty}\left(\dfrac{\tilde{C}_{n+1}}{\tilde{C}_n} - 1\right) \leq \sum_{n= 0}^{\infty}(\tilde{Q}_n-\tilde{Q}_{n+1})< +\infty,$$
where the latter sum is convergent as it is a telescoping sum with $\tilde{Q}_n \to 0$ as $n \to + \infty$. This assures that
$$\prod_{n = 0}^{\infty} \dfrac{\tilde{C}_{n+1}}{\tilde{C}_n} < \infty.$$

Summing up,
\begin{equation}
\label{eq:bound}
\dfrac{y_n}{x_{n+1}-x_n} \leq \dfrac{\tilde{y}_n}{\tilde{x}_{n+1}-\tilde{x}_n} \leq \tilde{C}_0\prod_{n = 0}^{+\infty} \dfrac{\tilde{C}_{n+1}}{\tilde{C}_n} =: C < +\infty,
\end{equation}
as desired.

\subsection{Proof of Theorem \ref{teo:perturb}}
\label{subsec:main}
We will use the ideas of the proof showed in Subsection \ref{subsec:secondproof} to prove Theorem \ref{teo:perturb}. In order to do so, write $f$ as
$$f(z) = z + \tilde{\beta} + p_1(z) + p_0(z), \quad z \in \H,$$
where
$$\omega(A) = \mu(A \cap (-\infty,0]), \quad \dfrac{d\nu}{d\mu}(t) = (1+t^2)\chi_{(0,+\infty)}(t),$$
and
$$\tilde{\beta} = \beta -\int_{(0,+\infty)}td\mu(t), \quad p_1(z) = \int_{\R}\dfrac{1+tz}{t-z}d\omega(t), \quad p_0(z) = \int_{\R}\dfrac{d\nu(t)}{t-z}.$$

Take $z_0 \in \H$. As usual, write $z_n = x_n + iy_n = f^n(z_0)$. The measure $\nu$ is finite by hypothesis, and then, if $y_0 \geq 1$, one has
$$\abs{\text{Re}(p_0(z_n))} \leq \int_{\R}\dfrac{d\nu(t)}{\abs{t-z_n}} \leq \dfrac{\nu(\R)}{y_n} \leq \nu(\R).$$ 
Thus, let $\beta^* \in \R$ be such that
\begin{equation}
\label{eq:betastar}
\beta^* < \tilde{\beta} - \nu(\R) \leq \tilde{\beta} - \abs{\text{Re}(p_0(z_n))}, \quad \text{for all } n \in \N.
\end{equation}
Define the function $\tilde{f} \colon \H \to \H$ given by
$$\tilde{f}(z) = z + \beta^* + p_1(z), \quad z \in \H.$$
Note that $\tilde{f}$ satisfies the hypothesis of Proposition \ref{prop:weak}, by construction. Let $b \in \R$ be as in Lemma \ref{lemma:properties} (applied to $\tilde{f}$). For any $x > b$, let us define
$$\tilde{P}(x) = \beta^* + \int_{(-\infty,0]}\dfrac{1+tx}{t-x}d\omega(t), \quad \tilde{Q}(x) = \int_{(-\infty,0]}\dfrac{1+t^2}{(t-x)^2}d\omega(t).$$
Again, given $\tilde{x}_0 > b$ and $\tilde{y}_0 > 0$ , define
$$\tilde{x}_{n+1} = \tilde{x}_n + \tilde{P}_n, \quad \tilde{y}_{n+1} = \tilde{y}_n(1+\tilde{Q}_n),$$
where $\tilde{P}_n = \tilde{P}(\tilde{x}_n)$ and $\tilde{Q}_n = \tilde{Q}(\tilde{x}_n)$. As shown in \eqref{eq:bound}, there exists $\tilde{C} > 0$ with
\begin{equation}
\label{eq:tildeC}
\dfrac{\tilde{y}_n}{\tilde{x}_{n+1}-\tilde{x}_n} < \tilde{C}.
\end{equation}

As in the previous subsection, the following relations among $x_n$, $y_n$, $\tilde{x}_n$ and $\tilde{y}_n$ hold:
\begin{lemma}
\label{lemma:main}
Suppose that $b < \tilde{x}_0 \leq x_0$ and $\tilde{y}_0 \geq y_0 \geq 1$. Then, for all $n \in \N$,
\begin{enumerate}[\hspace{0.5cm}\normalfont(a)]
\item $0 < \tilde{x}_{n+1}-\tilde{x}_n \leq x_{n+1}-x_n$, and $\tilde{x}_n \leq x_n$,
\item There exists $K > 0$ such that $K \tilde{y}_n \geq y_n$.
\end{enumerate}
\begin{proof}
(a) Due to the monotonicity properties showed in Lemma \ref{lemma:properties}, one has
$$\tilde{x}_{n+1}-\tilde{x}_n = \tilde{P}_n = \beta^* + p_1(\tilde{x}_n) \geq \beta^* + p_1(b) > 0.$$

Let us now proceed by induction. Suppose that $\tilde{x}_n \leq x_n$ and $\tilde{x}_{n+1}-\tilde{x}_n \leq x_{n+1}-x_n$. In that case,
$$\tilde{x}_{n+1} = \tilde{x}_n+(\tilde{x}_{n+1}-\tilde{x}_n) \leq x_n + (x_{n+1}-x_n) = x_{n+1}.$$
Furthermore, using the monotonicity and \eqref{eq:betastar},
\begin{align*}
\tilde{x}_{n+2} - \tilde{x}_{n+1} & = \beta^* + p_1(\tilde{x}_{n+1}) \leq \tilde{\beta} + p_1(\tilde{x}_{n+1}) + \text{Re}(p_0(z_{n+1})) \\
& \leq \tilde{\beta} + p_1(x_{n+1}) + \text{Re}(p_0(z_{n+1})) \\
& \leq \tilde{\beta} + \text{Re}(p_1(z_{n+1})) + \text{Re}(p_0(z_{n+1})) = x_{n+2}-x_{n+1}.
\end{align*}

(b) Define
$$Q_n = \dfrac{y_{n+1}-y_n}{y_n},$$
and note that
$$y_{n+1} = y_n(1+Q_n), \quad \tilde{y}_{n+1} = \tilde{y}_n(1+\tilde{Q}_n).$$
Then,
$$\dfrac{y_n}{\tilde{y}_n} = \dfrac{y_0}{\tilde{y}_0}\prod_{k = 0}^{n-1}\dfrac{1+Q_k}{1+\tilde{Q}_k} = \dfrac{y_0}{\tilde{y}_0}\prod_{k = 0}^{n-1}\left(1+\dfrac{Q_k-\tilde{Q}_k}{1+\tilde{Q}_k}\right) \leq \dfrac{y_0}{\tilde{y}_0}\prod_{k = 0}^{n-1}\left(1+\max\left\lbrace\dfrac{Q_n-\tilde{Q}_n}{1+\tilde{Q}_n},0\right\rbrace\right).$$

We claim that
$$\prod_{n = 0}^{\infty}\left(1+\max\left\lbrace\dfrac{Q_n-\tilde{Q}_n}{1+\tilde{Q}_n},0\right\rbrace\right) < \infty,$$
and so the proof follows. To prove the claim, it is enough to show that
$$\sum_{\substack{n = 0 \\ Q_n - \tilde{Q}_n \geq 0}}^{\infty}\dfrac{Q_n-\tilde{Q}_n}{1+\tilde{Q}_n} < +\infty.$$
To do so, note that
\begin{align*}
Q_n & = \int_{(-\infty,0]}\dfrac{1+t^2}{(t-x_n)^2+y_n^2}d\omega(t) + \int_{(0,+\infty)}\dfrac{d\nu(t)}{(t-x_n)^2+y_n^2} \\
& \leq \int_{(-\infty,0]}\dfrac{1+t^2}{(t-x_n)^2}d\omega(t) + \int_{(0,+\infty)}\dfrac{d\nu(t)}{(t-x_n)^2+y_0^2} \\
& \leq \int_{(-\infty,0]}\dfrac{1+t^2}{(t-\tilde{x}_n)^2}d\omega(t) + \int_{(0,+\infty)}\dfrac{d\nu(t)}{(t-x_n)^2+y_0^2} = \tilde{Q}_n + \int_{\R}\dfrac{d\nu(t)}{(t-x_n)^2+y_0^2}.
\end{align*}
Therefore,
$$\sum_{\substack{n = 0 \\ Q_n - \tilde{Q}_n \geq 0}}^{\infty}\dfrac{Q_n-\tilde{Q}_n}{1+\tilde{Q}_n} \leq \sum_{\substack{n = 0 \\ Q_n - \tilde{Q}_n \geq 0}}^{\infty}(Q_n-\tilde{Q}_n) \leq \sum_{n = 0}^{\infty}\int_{\R}\dfrac{d\nu(t)}{(t-x_n)^2+y_0^2},$$
and so it is enough to prove that
$$\sum_{n = 0}^{\infty}\int_{\R}\dfrac{d\nu(t)}{(t-x_n)^2+y_0^2} < +\infty.$$

On the one hand, note that the measure $\nu$ is finite, by hypothesis. On the other hand, one has
$$x_{n+1}-x_n \geq \tilde{x}_{n+1}-\tilde{x}_n = \beta^* + p_1(\tilde{x}_n) \geq \beta^* + p_1(b) > 0.$$
Thus, define $M = \beta^* + p_1(b) > 0$, and notice that $x_n \geq Mn + x_0$ and $x_n-x_m \geq M(n-m)$ for all $n,m \in \N$ with $n \geq m$. 

Assume that $0 < t \leq x_0$, and notice that
$$\sum_{n=0}^{\infty}\dfrac{1}{(t-x_n)^2+y_0^2} \leq \sum_{n=0}^{\infty}\dfrac{1}{(t-Mn-x_0)^2+y_0^2} \leq \sum_{n=0}^{\infty}\dfrac{1}{M^2n^2+y_0^2} < +\infty.$$
Otherwise, if $t > x_0$, it is possible to (uniquely) find $N \geq 0$ such that $t \in [x_N,x_{N+1})$. Then, for $0 \leq n \leq N$,
$$t-x_n = (t - x_N) + (x_N-x_n) \geq x_N-x_n \geq (N-n)M.$$
On the other hand, if $n > N$, then
$$x_n - t = (x_n-x_{N+1})+(x_{N+1}-t) \geq x_n-x_{N+1} \geq (n-N-1)M.$$
Joining both estimates,
\begin{align*}
\sum_{n=0}^{\infty}\dfrac{1}{(t-x_n)^2+y_0^2} & = \sum_{n = 0}^N\dfrac{1}{(t-x_n)^2+y_0^2} +\sum_{n = N+1}^{\infty}\dfrac{1}{(t-x_n)^2+y_0^2} \\
& \leq \sum_{n = 0}^N\dfrac{1}{(N-n)^2M^2+y_0^2} +\sum_{n = N+1}^{\infty}\dfrac{1}{(n-N-1)^2M^2+y_0^2} \\
& \leq 2\sum_{n=0}^{\infty}\dfrac{1}{n^2M^2+y_0^2} < +\infty.
\end{align*}

Combining all, we get
$$\sum_{n = 0}^{\infty}\int_{\R}\dfrac{d\nu(t)}{(t-x_n)^2+y_0^2} \leq 2\nu(\R)\sum_{n=0}^{\infty}\dfrac{1}{n^2M^2+y_0^2} < +\infty,$$
and we are done.
\end{proof}
\end{lemma}

We are now able to prove Theorem \ref{teo:perturb}. As before, the limit
$$L = \lim_{n \to + \infty}\dfrac{x_{n+1}-x_n}{y_n}$$
exists, as seen in Theorem \ref{teo:Pomm}. But, as seen in Lemma \ref{lemma:main}, $x_{n+1}-x_n > 0$, and so $L \geq 0$. Remember that $f$ is of PHS if and only if $L > 0$, that is, if and only if there exists some constant $C > 0$ such that
$$\dfrac{y_n}{x_{n+1}-x_n} \leq C.$$
But using \eqref{eq:tildeC},
$$\dfrac{y_n}{x_{n+1}-x_n} \leq K \dfrac{\tilde{y}_n}{\tilde{x}_{n+1}-\tilde{x}_n} \leq K\tilde{C} < +\infty.$$

\section{Examples}
\label{sec:examples}
Let $m$ denote the Lebesgue's measure on $\R$.
\begin{example}
Let $f \colon \H \to \H$ be the parabolic map given by
$$f(z) = z + \beta - \dfrac{1}{z+i}, \quad z \in \H,$$
where $\beta \in \R$. Define $\mu$ as the positive finite measure on $\R$, which is absolutely continuous with respect to $m$, and given by
$$\dfrac{d\mu}{dm}(t) = \dfrac{1}{\pi}\dfrac{1}{(1+t^2)^2}, \quad t \in \R.$$
Note that
$$\int_{\R}t^2 d\mu(t) < \infty, \quad \int_{\R}td\mu(t) = 0.$$
Also,
$$\int_{\R}\dfrac{1+tz}{t-z}d\mu(t) = -\dfrac{1}{z+i}, \quad z \in \H.$$
By Theorem \ref{teo:t1end}, $f$ is of 0HS if and only if $\beta = 0$.

Also note that
$$\angle\lim_{z\to\infty}z(f(z)-z) =
\begin{cases}
\infty, & \beta \neq 0, \\
-1, & \beta = 0,
\end{cases}$$
as deduced in Proposition \ref{prop:t2charac}.
\end{example}
\begin{example}
Let $f \colon \H \to \H$ be the parabolic map given by
$$f(z) = z + \beta -\dfrac{\log(1-z)}{z}-\dfrac{\pi}{4}, \quad z \in \H,$$
where we are using the principal argument to define the logarithm and $\beta \in \R$. Define $\mu$ as the positive finite measure on $\R$, which is absolutely continuous with respect to $m$, and given by
$$\dfrac{d\mu}{dm}(t) = \dfrac{1}{1+t^2}\dfrac{1}{t}\chi_{[1,+\infty)}(t), \quad t \in \R.$$
Note that
$$\int_{\R}\abs{t} d\mu(t) < \infty, \quad \int_{\R} td\mu(t) = \dfrac{\pi}{4},$$
$$\int_{(-\infty,0)}t^2d\mu(t) = 0, \quad \int_{(0,+\infty)}t^2d\mu(t) = \infty.$$
Also,
$$\int_{\R}\dfrac{1+tz}{t-z}d\mu(t) = -\dfrac{\log(1-z)}{z}-\dfrac{\pi}{4}, \quad z \in \H.$$
By Theorem \ref{teo:t1end}, $f$ is of 0HS if and only if $\beta \geq \pi/4$.
\end{example}
\begin{example}
Let $f \colon \H \to \H$ be a parabolic map given by
$$f(z) = z + \beta + \tan(z), \quad z \in \H,$$
where $\beta \in \R$. Define $\mu$ as the positive finite measure on $\R$, which is singular with respect to $m$, and given by
$$\mu = \sum_{n \in \Z}\dfrac{\delta_{p_n}}{1+p_n^2},$$
where
$$p_n = \dfrac{\pi}{2}(2n+1), \quad n \in \Z.$$
Note that $\mu$ is symmetric, and
$$\int_{\R}\abs{t} d\mu(t) = \infty.$$
Also,
\begin{align*}
\int_{\R}\dfrac{1+tz}{t-z}d\mu(t) & = \int_{(0,+\infty)}\dfrac{2z(1+t^2)}{t^2-z^2}d\mu(t) = \sum_{n = 0}^{\infty} \dfrac{2z}{p_n^2-z^2} \\
& = \sum_{n= 0}^{\infty}\dfrac{8z}{(2n+1)^2\pi^2-4z^2} = \tan(z), \quad z \in \H.
\end{align*}
This is a known representation of the tangent function as a series (see \cite[Formula (3), p. 329]{TangentSerieBook}). By Theorem \ref{teo:symmetric_conc} or Theorem \ref{teo:symmetric}, $f$ is of 0HS for every $\beta \in \R$.
\end{example}
\begin{example}
Let $f \colon \H \to \H$ be a parabolic map given by
$$f(z) = z + \beta +\log(z), \quad z \in \H,$$
where we are using the principal argument to define the logarithm and $\beta \in \R$. Define $\mu$ as the positive finite measure on $\R$, which is absolutely continuous with respect to $m$, and given by
$$\dfrac{d\mu}{dm}(t) = \dfrac{1}{1+t^2}\chi_{(-\infty,0)}(t), \quad t \in \R.$$
Note that
$$\int_{\R}\abs{t} d\mu(t) = \infty, \quad \mu((0,+\infty)) = 0.$$
Also,
$$\int_{\R}\dfrac{1+tz}{t-z}d\mu(t) = \log(z), \quad z \in \H.$$
By Theorem \ref{teo:perturb}, $f$ is of PHS for every $\beta \in \R$.
\end{example}

\end{document}